\UseRawInputEncoding
\documentclass{amsart}
\usepackage{sansmath}
\usepackage{amssymb,amsthm,amsmath}
\usepackage[numbers,sort&compress]{natbib}
\usepackage{color}
\usepackage{graphicx}
\usepackage{stmaryrd}
\usepackage{tikz}
\usepackage{caption}
\usepackage[latin,english]{babel}
\usepackage{braket}
\usepackage[toc,page]{appendix}

\title[ ]{The stability of Sobolev norms for the linear wave equation with unbounded perturbations}

\author{Yingte Sun}
\address[Y. Sun]{School of Mathematical Sciences,
Yangzhou University,
yangzhou,
China} \email{sunyt15@fudan.edu.cn}

\DeclareSymbolFont{largesymbol}{OMX}{yhex}{m}{n}
\DeclareMathAccent{\Widehat}{\mathord}{largesymbol}{"62}

\newcommand{\R}{\mathbb{R}}

\newcommand{\T}{\mathbb{T}}

\theoremstyle{plain}
\newtheorem{thm}{Theorem}[section]
 \newtheorem{cor}[thm]{Corollary}
 \newtheorem{lem}[thm]{Lemma}
 \newtheorem{prop}[thm]{Proposition}
 \theoremstyle{definition}
 \newtheorem{defn}[thm]{Definition}
 \theoremstyle{remark}
 \newtheorem{rem}[thm]{Remark}
 
 \numberwithin{equation}{section}

\begin{document}


\begin{abstract}
In this paper, we prove that  the Sobolev norm of solutions of the linear wave equation with unbounded perturbations of order one stay bounded for the all time. The main proof is based on the KAM reducibility of the linear wave equation. To the best of our knowledge, this is the first reducibility result of the linear wave equation with general quasi-periodic unbounded potentials on the torus.
\end{abstract}

\maketitle
\tableofcontents
\section{Introduction }
In this paper, we consider a linear wave equation with unbounded quasi-periodic perturbations of the form
\begin{equation}\label{wave1}
\partial_{tt}u-\partial_{xx}u+\mathrm{m}u+\mathcal{W}(\omega t)u=0, \quad t \in\mathbb{R}, \quad  x\in \mathbb{T}=\mathbb{R}/2\pi \mathbb{Z},
\end{equation}
where $\mathcal{W}(\omega t)$ is a pseudo-differential operator of order one, and quasi-periodic in time with frequencies $\omega \in \mathcal{O}:=[1,2]^d$. The mass $\mathrm{m}$ is positive. We prove that the Sobolv norms of solutions $(u,u_{t})$ of the equation \eqref{wave1} are uniformly bounded for a large subset of $\mathcal{O}$. The main proof is based on a quantitative reducibility result of the wave equation, in which we construct a bounded and time quasi-periodic transformation on space $\mathcal{H}^r_x \times \mathcal{H}^r_x$ such that the original equation \eqref{wave1} can be transformed into a block diagonal and time independent one.

The problem of estimating the high Sobolev norm of linear partial differential equations  has been widely studied. The two remarkable results were obtained by  Bourgain \cite{Bou99,Bou99.1} for the free Schr\"odinger equation with
time dependent potential
\begin{equation}
\mathrm{i}\partial_t u=-\Delta u+V(t,x)u,
\end{equation}
on the $d$-dimensional torus. Bougain \cite{Bou99} derived a $\langle t\rangle^{\epsilon}$ upper bound of the Sobolev norm of solutions for smooth and bounded time dependent potentials. When the potential $V$ is analytic and time quasi-periodic, Bourgain \cite{Bou99.1} proved that the Sobolev norm of solution grows like a power of $log(t)$. The result obtained in \cite{Bou99} has  been extended by Delort \cite{Delort2010} and Berti-Maspero\ cite{berti19} to the Zoll manifolds and flat tours. The logarithmic bounds on Sobolev norms in \cite{Bou99.1} has been extended by Wang \cite{Wang08} to an analytic and bounded time dependent potential on $\mathbb{T}$.

However, Bourgain's original method can only deal with the bounded perturbations, especially the multiplicative potential $V(t,x)$. The first result on the Schr\"odinger-type equation with unbounded, time dependent perturbation of the form
 \begin{equation}
\mathrm{i}\partial_tu(t)=\mathrm{H}u(t)+\mathrm{P}(t) u,
\end{equation}
 is due to Maspero-Robert \cite{Mas2017}. The method in \cite{Mas2017} can be applied to the free Schr\"odinger  equation on Zoll manifolds with time smooth perturbations of order $m<1$, which provided a $\langle t\rangle^{\epsilon}$ upper bound of Sobolev norms of solutions. Based on the pseudo-differential operator technique, Bambusi-Gr\'ebert- Maspero-Robert \cite{Bam21} extended their results to more Schr\"odinger-type equations, including the free Schr\"odinger on Zoll manifolds with perturbations of order $m<2$.  In the meantime Montolto \cite{Mon2018.1,Mo20191} has independently studied the maximum order of perturbations for the Schr\"odinge-type equation on $\mathbb{T}$.  It's worth mentioning that, based on a  delicate  Quantum version Nehoroshev  theorem \cite{Bam22}, Bambusio-Langella-Montalto \cite{Bam20} proved a $\langle t \rangle^\varepsilon$ upper bound of Sobolev norm of solutions  for the free Schr\"odinger equation on the flat torus with unbounded perturbations of order $m<2$.

For the Schr\"odinger-type equation  with small  time quasi-periodic perturbations of the form
$$\mathrm{i}\partial_tu(t)=\mathrm{H}u(t)+\varepsilon \mathrm{P}(\omega t) u,$$
the KAM reducibility is a powerful tool to investigate the uniformly boundedness of the solutions in Sobolev space. For the bounded perturbations, we mention the results
of Eliasson-Kuksin \cite{L.H 09} which proved the reducibility of the Schr\"odinger equation on $\mathbb{T}^d$  and Gr\'ebert et al. \cite{B12,B24} which proved the
reducibility of the quantum harmonic oscillator on $\mathbb{R}^d$.The reducibility results imply that the Sobolev norms of solutions for the equation considered is uniformly bounded. For the unbounded perturbations, there are several papers devoting to the reducibility of some Schr\"odinger equations, such as the quantum harmonic oscillator \cite{Bam18,Bam171,Bam19,Bam018}, duffing oscillator \cite{Liu09}, relativistic Schr\"odinger equation on torus  \cite{S21} and the free schr\"odinger equation on Zoll manifolds \cite{F.G2019,F.G2020}.

Compared with the enormous results of Schr\"odinger-type equations, there are few results concerning the growth of Sobolev nrom of solutions for the wave equation.  Estimating the high Sobolev norm of solutions for linear wave equations on compact manifolds is  much more subtle than  Schr\"odinger-type equations. Fang-Han-Wang \cite{Fang14} had constructed a small time periodic potential for the wave equation on the torus such that the Sobolev norms of solutions is bounded for all time. While, Bourgain \cite{Bou99.1} constructed a time periodic potential for the  wave equation,
which provoke exponential growth of Sobolev norm. In order to avoid such terrible upper bound, people has to  pay more attention to the wave equation with time quasi-periodic perturbation. Naturally, the KAM reducibility  becomes one of the main research methods. For the bounded perturbations, we mention the results of Li\cite{Li2020} and Liang \cite{Liang2018} which proved the reducibility of  wave equation on the torus with small time quasi-periodic  multiplicative potential. Maspero \cite{L.M2019} proved the reducibility of the wave equation with non-small time quasi-periodic  multiplicative potential.  For the unbounded perturbations,  Montolto \cite{Mo2019}, Sun et al. \cite{S19} studied the wave equation with some special unbounded perturbations.

It should be mentioned that the unbounded perturbations consider in \cite{Mo2019,S19} is in the special form of  $V(\omega t)\Delta$, which can be obtained by linearizing some nonlinear equations \cite{Mo20191}. People are more concerned with general unbounded perturbations. From the viewpoint of KAM theory, if the order of perturbations is strictly smaller than one, the KAM reducibility for such  wave equation is straightforward. If the order of perturbations is equal to one, some serious problems are occurs in the measure estimate in KAM iteration. The similar problems are resolved by Berti-Biasco-Procesi \cite{Berti13}, in which they obtained some quasi-periodic solutions of the Hamiltonian derivative wave equation. In order to estimate the number of non-resonance conditions, they introduced the ``quasi-T\"oplitz" property of the perturbations to get the higher order asymptotic decay estimate of normal frequencies.
However, the property of  momentum conservation of the equation is indispensable for  preserveing the ``quasi-T\"oplitz" property in KAM iteration. The property is missing from the wave equation \eqref{wave1} in the present paper. Therefore, this paper adopts a completely different method, that is the method of pseudo-differential operator.

For the  Schr\"odinger-type equation with unbounded perturbations, the method of  pseudo-differential operator can effectively smoothing the perturbations, so as to avoid a series of difficulties caused by the order of perturbations. However, such skills are almost useless for the wave equation. One of the main reasons is that  the wave equation with Hamiltonian form can be seen as a $2\times 2$ matrix valued, Schr\"odinger-type equation \eqref{2}. For the matrix-valued pseudo-differential operator, the commutator of two  matrix-valued pseudo-differential operators can not gains one derivative. Therefore, we can not transform the whole perturbation in equation \eqref{2} into a smoothing one. The main novelty of the present paper is that we find a delicate bounded transformation such that the original perturbations $\mathbf{K}(\omega t)$ in  equation \eqref{2} can be transformed  into a new one $\mathbf{P}(\omega t)$ in  equation \eqref{wave3}, where the diagonal part are smoothing operators and the off-diagonal part are still  bounded operators.
Furthermore, such structure can be maintained in the KAM iteration. Under these conditions we can also get the higher order asymptotic decay estimate of eigenvalues.

\begin{rem} If we have a good control of the matrix decay norm of the operator $e^{\mathrm{i} \mathcal{G}}$, where $\mathcal{G}$ is a self-adjoint, pseudo-differential operator of order $0< m < \frac{1}{2}$, the method presented in this paper may be extended to the case that the order of perturbations  is less than $3/2$. Form the Lemma 3.4 in \cite{Bam20}, we known that the operator $e^{\mathrm{i}  \mathcal{G}}$ is  bounded in Sobolev space $\mathcal{H}^r_{x}$. But the information about its  matrix decay norm is missing, which is essential to the KAM iteration.
\end{rem}

\begin{rem}
In addition, the reducibility problem of wave equation \eqref{wave1} in this paper is the  cornerstone  of some further works. Inspired by \cite{bal19,Baldi2}, we can use the quantitative reducibility result in this paper to explore the existence of Sobolev, linearly stable, quasi-periodic solutions of the following derivatives wave equations.

\textbf{$\bullet$ The Hamiltonian derivatives wave equations with quasi-periodic force}
\begin{equation}
\partial_{tt}u-\partial_{xx}u+\mathrm{m}u+f(\omega t,x,Du)=0,\quad D=\sqrt{-\partial_{xx}+\mathrm{m}}, \quad  x\in \mathbb{T}.
\end{equation}

\textbf{$\bullet$ The autonomous derivatives wave equations}
\begin{equation}
\partial_{tt}u-\partial_{xx}u+\mathrm{m}u+a(x)f(Du)=0,\quad D=\sqrt{-\partial_{xx}+\mathrm{m}},\quad  x\in \mathbb{T}.
\end{equation}
\end{rem}

The paper is organized as follows: In section 2, we introduce some important  definitions of pseudo-differential operator, so that we can precisely state our main results. In section 3, we introduce some norms of infinite dimensional
matrix, such that the KAM iterations in section 5 can be well understood. In section 4, we
introduce the the symbolic calculus of pseudo-differential operators in \cite{Mo20191}, such that the
diagonal part of the  perturbation $\mathbf{K}(\omega t)$ in equation \eqref{2}  can be reduced to a operator of order $-1$. In section 5, we give
a block-diagonal reducibility result for the equation \eqref{2}.  In Section 6, we
conclude the proof of Theorem \ref{main1} and  the Corollary \ref{main2}.

Notations: In the present paper, we denote the notation $A\lesssim_s  B$ as $A \leq C(s)B$, where $C(s)$ depends on the data of the problem, namely  the Sobolev index $s$, the number $d$ of time frequencies, the diophantine exponent $\tau>0$  in the non-resonance conditions, which will be required along the proof.
\section{Main result}
Given a Function $f:\mathcal{O}\mapsto E: \omega \mapsto f(\omega)$, where $(E, \|\cdot\|_E)$ is Banach space and $\omega \in \mathcal{O}$. We define the sup-norm and  lipschitz semi-norm as
\begin{equation}
\|f\|^{\mathrm{sup}}_{E,\mathcal{O}}:=\sup_{\omega \in \mathcal{O}}\|f(\omega)\|_{E}, \quad  \|f\|^{lip}_{E,\mathcal{O}}:=\sup_{\stackrel{\omega_1,\omega_2 \in \mathcal{O}}{\omega_1 \neq \omega_2}}\frac{\|f(\omega_1)-f(\omega_2)\|_{E}}{|\omega_1-\omega_2|}.
\end{equation}
For any $\gamma >0$, we define the Lipschitz-norm
\begin{equation}
\|f\|^{\gamma}_{E,\mathcal{O}}:=\|f\|^{\mathrm{sup}}_{E,\mathcal{O}}+\gamma \|f\|^{lip}_{E,\mathcal{O}}.
\end{equation}
For notation convenience, we omit to write the set $\mathcal{O}$.
\subsection{Function space and pseudo-differential operators}\

\textbf{Sobolev space:}

For any function $u(x) \in L^2(\mathbb{T})$, it can be written as
$$u(x)=\sum_{j\in \mathbb{Z}} \hat{u}_je ^{\mathrm{i}j\cdot x}, \quad  \hat{u}_j=\frac{1}{2\pi}\int_{\mathbb{T}}u(x)e^{-\mathrm{i}j\cdot x}dx.$$
The Sobolev space  $\mathcal{H}^{s}_x$ is defined by
$$\mathcal{H}^{s}_x:=\Big\{u(x)\in L^2(\mathbb{T}):\Big| \|u(x)\|^{2}_{\mathcal{H}^{s}_x}=\sum_{j \in \mathbb{Z}}\langle j \rangle^{2s}|\hat{u}_j|^2 < +\infty \Big\},$$
where $\langle j\rangle=\max\{1,|j|\}$.

For any functions $u(\theta,x) \in L^2(\mathbb{T}^d \times \mathbb{T})$, it can be regarded as a $\theta -$dependent family of functions $u(\theta,\cdot) \in L^2(\mathbb{T})$. We can expand in Fourier series as
$$u(\theta,x)=\sum_{j\in \mathbb{Z}}\hat{u}_j(\theta)e^{\mathrm{i}j \cdot x}=\sum_{(\ell,j)\in \mathbb{Z}^{d+1}}\hat{u}_{j}(\ell)e^{\mathrm{i}(j\cdot x+\ell \cdot x)},$$
where
$$\hat{u}_j(\theta)=\frac{1}{2\pi} \int_{\mathbb{T}}u(\theta,x)e^{-\mathrm{i}j\cdot x}, \quad \hat{u}_j(\ell)=\frac{1}{(2\pi)^{d+1}}\int_{\mathbb{T}^{d+1}}u(\theta,x)e^{-\mathrm{i}(j\cdot x+\ell\cdot \theta)}dxd\theta.$$
The Sobolev space $\mathcal{H}^s(\mathbb{T}^d\times \mathbb{T})$ is defined by
$$\mathcal{H}^{s}(\mathbb{T}^{d}\times \mathbb{T}):=\Big\{u \in L^2(\mathbb{T}^d\times \mathbb{T}):\Big| \|u\|^2_s= \sum_{(\ell,j )\in \mathbb{Z}^{d+1}} \langle \ell ,j \rangle^{2s}|\hat{u}_j(\ell)|^2 < +\infty  \Big\},$$
where $\langle\ell,  j\rangle=\max\{1,|j|,|\ell|\}$.

Notation: In the rest of the paper, we fix
$$s_0:=\big[\frac{d+1}{2}\big]+1,$$
where for any real number $x\in\mathbb{R}$, we denote by $[x]$ its integer part.

\textbf{Pseudo-differential operators:}
\begin{defn}\label{defnpd}(Pseudo-differential operators and symbols)
Let $m \in \mathbb{R}, s\geq s_0, \alpha \in \mathbb{N}$, we say that an operator $\mathcal{A}=\mathcal{A}(\theta)$ is in the class $\mathcal{OPS}^m_{s,\alpha}$, if there exists a function $a=a(\theta,x,\xi):\mathbb{T}^d\times \mathbb{T}\times \mathbb{R}\mapsto \mathbb{C}$, differentiable $\beta$ times in the variables $\xi$, such that
$$\mathcal{A}u(x)=\mathrm{Op}(a)u(x)=\sum_{\xi \in \mathbb{Z}}a(\theta,x,\xi)\hat{u}(\xi)e^{\mathrm{i}x\cdot\xi},\quad \forall u\in \mathcal{H}^{0}_{x},$$
and
$$|\mathcal{A}|_{m,s,\alpha}:=\sup_{|\beta|\leq\alpha}\sup_{\xi \in \R}\|\partial^{\beta}_{\xi}a(\theta,x,\xi)\|_s\langle \xi\rangle^{-m+\beta}.$$
In that case, we say that $a(\theta,x,\xi)$ is in the class $\mathcal{S}^m_{s,\alpha}$. The operator $\mathcal{A}$ is said to be a pseudo-differential operator of order $m$, and the function $a$ is symbol.

If $\mathcal{A}:=\mathcal{A}(\lambda)$ is depending in a Lipschitz way on some parameter $\omega \in \mathcal{O} \subseteq \mathbb{R}^d$, we set
\begin{equation}
|\mathcal{A}|^{\gamma}_{m,s,\alpha}=|\mathcal{A}|^{\gamma,\mathcal{O}}_{m,s,\alpha}:=\sup_{\omega\in \mathcal{O}}|\mathcal{A}|_{m,s,\alpha}+\gamma \sup_{\omega_1, \omega_2\in \mathcal{O}}\frac{|\mathcal{A}(\omega_1)-\mathcal{A}(\omega_2)|_{m,s,\alpha}}{|\omega_1-\omega_2|}.
\end{equation}
\end{defn}

\begin{lem}\label{estpd}(Lemmata 2.13, 2.15 in \cite{bert20})Let $s\geq s_0$, $m,m' \in \R$, $\alpha \in \mathbb{N}$.

$\mathbf{1:}$ Let $\mathcal{A}:=\mathrm{Op}(a)\in \mathcal{OPS}^{m}_{s+|m|+\alpha,\alpha}$,$\mathcal{B}:=\mathrm{Op}(b)\in \mathcal{OPS}^{m'}_{s,\alpha}$, then the composition $\mathcal{AB}$ belongs to $\mathcal{OPS}^{m+m'}_{s,\alpha}$,and
\begin{equation}
|\mathcal{AB}|^{\gamma}_{m+m',s,\alpha}\lesssim _{s,m,\alpha} |\mathcal{A}|^{\gamma}_{m,s,\alpha}|\mathcal{B}|^{\gamma}_{m',s_0+|m|+\alpha,\alpha}+|\mathcal{A}|^{\gamma}_{m,s_0,\alpha}|\mathcal{B}|^{\gamma}_{m',s+|m|+\alpha,\alpha}.
\end{equation}

$\mathbf{2:}$ Let $\mathcal{A}:=\mathrm{Op}(a)\in \mathcal{OPS}^{m}_{s,\alpha+1}$, $\mathcal{B}:=\mathrm{Op}(b)\in \mathcal{OPS}^{m'}_{s+|m|+\alpha+2,\alpha}$. Then
$$\mathcal{AB}=\mathrm{Op}(a(\theta,x,\xi)b(\theta,x,\xi))+\mathfrak{R}_{\mathcal{AB}}, \quad \mathfrak{R}_{\mathcal{AB}}\in \mathcal{OPS}^{m+m'+1}_{s,\alpha},$$
where the Reminder $\mathfrak{R}_{AB}$ satisfies
$$|\mathfrak{R}_{\mathcal{AB}}|^{\gamma}_{m+m'-1,s,\alpha} \lesssim_{s,m,\alpha}|\mathcal{A}|^{\gamma}_{m,s,\alpha+1}|\mathcal{B}|^{\gamma}_{m',s_0+|m|+2,\alpha}
+|\mathcal{A}|^{\gamma}_{m,s_0,\alpha+1}|\mathcal{B}|^{\gamma}_{m',s+|m|+2,\alpha}.$$

$\mathbf{3:}$ Let $\mathcal{A}:=\mathrm{Op}(a)\in \mathcal{OPS}^{m}_{s,\alpha+2}$, $\mathcal{B}:=\mathrm{Op}(b)\in \mathcal{OPS}^{m}_{s+|m|+\alpha+4,\alpha}$. Then
$$\mathcal{AB}=\mathrm{Op}(a(\theta,x,\xi)b(\theta,x,\xi)-\mathrm{i}\partial_{\xi}a(\theta,x,\xi)\partial_{x}b(
\theta,x,\xi))+\mathfrak{R}_{2,\mathcal{AB}}, \quad \mathfrak{R}_{2,\mathcal{AB}}\in \mathcal{OPS}^{m+m'-2}_{s,\alpha},$$
where the Reminder $\mathfrak{R}_{2,\mathcal{AB}}$ satisfies
$$|\mathfrak{R}_{2,\mathcal{AB}}|^{\gamma}_{m+m'-2,s,\alpha} \lesssim_{s,m,\alpha}|\mathcal{A}|^{\gamma}_{m,s,\alpha+2}|\mathcal{B}|^{\gamma}_{m',s_0+|m|+4,\alpha}
+|\mathcal{A}|^{\gamma}_{m,s_0,\alpha+2}|\mathcal{B}|^{\gamma}_{m',s+|m|+4,\alpha}.$$
\end{lem}
\begin{rem}\label{pse} \

$\bullet$ From item $\mathbf{2}$ in Lemma \ref{estpd}, if $\mathcal{A}:=\mathrm{Op}(a)\in \mathcal{OPS}^{m}_{s+|m'|+\alpha+2,\alpha+1}$, $\mathcal{B}:=\mathrm{Op}(b)\in \mathcal{OPS}^{m'}_{s+|m|+\alpha+2,\alpha+1}$, then, the commutator $[\mathcal{A},\mathcal{B}]:=\mathcal{A}\mathcal{B}-\mathcal{B}\mathcal{A} \in \mathcal{OPS}^{m+m'-1}_{s,\alpha}$, and
\begin{equation}
\begin{split}
|[\mathcal{A},\mathcal{B}]|^{\gamma}_{m+m'-1,s,\alpha}\lesssim_{m,m',\alpha,s}&|\mathcal{A}|^{\gamma}_{m,s+|m'|+\alpha+2,\alpha+1}|\mathcal{B}|^{\gamma}_{m,s_0+|m|+\alpha+2,\alpha+1}\\
&+|\mathcal{A}|^{\gamma}_{m,s_0+|m'|+\alpha+2,\alpha+1}|\mathcal{B}|^{\gamma}_{m,s+|m|+\alpha+2,\alpha+1}.
\end{split}
\end{equation}

$\bullet$ From item $\mathbf{3}$ in Lemma \ref{estpd}, if $\mathcal{A}:=\mathrm{Op}(a)\in \mathcal{OPS}^{m}_{s+|m'|+\alpha+4,\alpha+2}$, $\mathcal{B}:=\mathrm{Op}(b)\in \mathcal{OPS}^{m'}_{s+|m|+\alpha+4,\alpha+2}$, then, the commutator $[\mathcal{A},\mathcal{B}]:=\mathrm{Op}(-\mathrm{i}\{a,b\}+r_{a,b})$, where $\{a,b\}=\partial_{\xi}a\partial_{x}b-\partial_{x}a\partial_{\xi}b$ and $\mathrm{Op}(r_{a,b})\in \mathcal{OPS}^{m+m'-2}_{s,\alpha}$ satisfies

\begin{equation}
\begin{split}
|\mathrm{Op}(r_{a,b})|^{\gamma}_{m+m'-2,s,\alpha}\lesssim_{m,m',\alpha,s}&|\mathcal{A}|^{\gamma}_{m,s+|m'|+\alpha+4,\alpha+2}|\mathcal{B}|^{\gamma}_{m,s_0+|m|+\alpha+4,\alpha+2}\\
&+|\mathcal{A}|^{\gamma}_{m,s_0+|m'|+\alpha+4,\alpha+2}|\mathcal{B}|^{\gamma}_{m,s+|m|+\alpha+4,\alpha+2}.
\end{split}
\end{equation}

\end{rem}
\textbf{Adjoint of pseudo-differential operator}

Considering a $\theta$-dependent families of pseudo-differential operator $\mathcal{A}(\theta)$, the symbol of the adjoint operator $\mathcal{A}^*=\mathrm{Op}(a^*(\theta,x,\xi))$ is
\begin{equation}
a^*(\theta,x,\xi)=\overline{\sum_{j \in \mathbb{Z}}\hat{a}(\theta,j,\xi-j)e^{\mathrm{i}j\cdot x}}=\overline{\sum_{\ell \in \mathbb{Z}^d,j\in \mathbb{Z}}\hat{a}(\ell,j,\xi-j)e^{\mathrm{i}(j \cdot x+\ell \cdot \theta)}}
\end{equation}

\begin{lem}(Lemma 2.16 in \cite{bert20})
Let $\mathcal{A}=\mathrm{Op}(a) \in \mathcal{OPS}^{m}_{s+s_0+|m|,0}$, and dependent on the parameters $\omega \in  \mathcal{O}$. Then, the adjoint operator $\mathcal{A}^*$ satifies
\begin{equation}
|\mathcal{A}^*|^{\gamma}_{m,s,0} \lesssim_{m,s} |\mathcal{A}|^{\gamma}_{m,s+s_0+|m|,0}.
\end{equation}
\end{lem}
\begin{lem}\label{adjsym1}
Let $\mathcal{A}=\mathrm{Op}(a(\theta,x,\xi)) \in \mathcal{OPS}^{m}_{s+|m'|+2,\alpha}$ be a self-adjoint operator and $\mathcal{G}=\mathrm{Op}(g(\xi))$ be a real Fourier multiplies of order $m'$(independent of parameters $\omega$). We define a new operator $\mathcal{B}=\mathrm{Op}(g(\xi)\cdot a(\theta,x,\xi))$, that
\begin{equation}\label{self1}
|\mathcal{B}|^{\gamma}_{m+m',s+|m'|+2,\alpha}\lesssim_{s,\alpha}|\mathcal{A}|^{\gamma}_{m,s+|m'|+2,\alpha}.
\end{equation}
Also, the operator $\mathcal{B}-\mathcal{B}^* \in \mathcal{OPS}^{m+m'-1}_{s,\alpha}$ and satisfies
\begin{equation}\label{self2}
|\mathcal{B}-\mathcal{B}^*|^{\gamma}_{m+m'-1,s,\alpha} \lesssim_{m,s,\alpha}|\mathcal{A}|^{\gamma}_{m,s+|m'|+2,\alpha}.
\end{equation}
\end{lem}

\begin{proof}
The estimate \eqref{self1} is a direct corollary of Definition \ref{defnpd}.

For the composition operator $\mathcal{B}$, one sees
$$\mathcal{B}=\mathrm{Op}(g(\xi)\cdot a(\theta,x,\xi))=\mathrm{Op}( a(\theta,x,\xi))\circ \mathrm{Op}(g(\xi))$$
and
$$\mathcal{B}^*=\mathrm{Op}^*(g(\xi)) \circ \mathrm{Op}^*( a(\theta,x,\xi)). $$
Since the operator $\mathcal{A}$ is self-adjoint and $g(\xi)$ is real, one has
\begin{equation}
\mathcal{B}^*=\mathrm{Op}(g(\xi)) \circ \mathrm{Op}( a(\theta,x,\xi)).
\end{equation}
From Lemma \ref{estpd},one gets that
$$\mathcal{B}^*=\mathrm{Op}(g(\xi)) \circ \mathrm{Op}( a(\theta,x,\xi))=\mathrm{Op}(g(\xi)\cdot a(\theta,x,\xi))+\mathcal{R}_{\mathcal{G,A}},$$
and
\begin{equation}
\begin{split}
|\mathcal{B}^*-\mathcal{B}|^{\gamma}_{m+m'-1,s,\alpha}
&=|\mathfrak{R}_{\mathcal{G},\mathcal{A}}|^{\gamma}_{m+m'-1,s,\alpha}\\
& \lesssim_{s,m',\alpha}|\mathcal{G}|^{\gamma}_{m',s,\alpha+1}|\mathcal{A}|^{\gamma}_{m,s_0+|m|+2,\alpha}
+|\mathcal{G}|^{\gamma}_{m',s_0,\alpha+1}|\mathcal{A}|^{\gamma}_{m,s+|m'|+2,\alpha} \\
& \lesssim_{s,m',\alpha}|\mathcal{A}|^{\gamma}_{m,s+|m'|+2,\alpha}.
\end{split}
\end{equation}
\end{proof}
For any symbol $a\in \mathcal{S}^m_{s,\alpha}$, we defined the average symbol $\langle a\rangle_{\theta,x}$ by
\begin{equation}
\langle a\rangle_{\theta,x}=\frac{1}{(2\pi)^{d+1}}\int_{\mathbb{T}^{d+1}}a(\theta,x,\xi)d\theta dx.
\end{equation}
Given $\omega \in \mathbb{R}^d$ and satisfies the non-resonance condition
\begin{equation}\label{non}
|\omega \cdot \ell \pm j|\geq \frac{\gamma}{\langle \ell \rangle^{\tau}},\quad  , \quad \forall(\ell, j) \in \mathbb{Z}^{d+1} \backslash \{0\}.
\end{equation}
We define the operator $(\omega \cdot \partial_{\theta} \pm \partial_{x})^{-1}$ by setting
$$(\omega \cdot \partial_{\theta}\pm \partial_x)^{-1}[1]=0, \quad  (\omega \cdot \partial_{\theta}\pm \partial_x)^{-1}(e^{\mathrm{i}(\ell \cdot\theta+j\cdot x)})=\frac{e^{\mathrm{i}(\ell \cdot\theta+j\cdot x)}}{\mathrm{i}(\omega\cdot \ell \pm j)},  \quad \forall(\ell, j) \in \mathbb{Z}^{d+1} \backslash \{0\}. $$

\begin{lem}(Lemma 2.8 in \cite{Mo20191})\label{adjsym2}Given a symbol $a \in \mathcal{S}^{m}_{s,\alpha}$,

\textbf{1}: $\langle a ^*\rangle_{\theta,x}=\overline{\langle a \rangle_{\theta,x}} =(\langle a\rangle_{\theta,x})^*.$

\textbf{2}: if $\omega$ satisfies the non-resonance condition \eqref{non}, then
$$(\omega\cdot \partial_{\theta} \pm \partial_x)^{-1} a^{*}=\big((\omega\cdot \partial_{\theta} \pm \partial_x)^{-1} a\big)^{*}.$$
\end{lem}
We define the operator $\sqrt{-\Delta}$ as follows, let  $\chi \in \mathcal{C}^{\infty}(\mathbb{R},\mathbb{R})$ be  a cut-off function satisfies
$$\chi(\xi):=\left\{
\begin{aligned}
&1  \quad  if \ |\xi|\geq 1, \\
&0 \quad  if \ |\xi| \leq \frac{1}{2}.
\end{aligned}
\right.
$$
We then define the operator $\sqrt{\Delta}$ as $\mathrm{Op}(|\xi|\chi(\xi)).$

\subsection{Main result}
Consider the perturbation $\mathcal{W}(\omega t)$ in equation \eqref{wave1}, we assume that

\textbf{Condition I}: $\mathcal{W}(\omega t)$ is a real, and self-adjoint linear operator.

\textbf{Condition II}: Set the symbol of the pseudo-differential operator $\mathcal{W}(\omega t)$ as $w(\theta,x,\xi)$, it satisfies
\begin{equation}
\langle w\rangle_{\theta,x}=\int_{\mathbb{T}^{d+1}}w(\theta,x,\xi)dx d\theta=a(\xi)\langle \xi \rangle+b(\xi),
\end{equation}
where $a(\xi)\in \Gamma=\{a_1^*,\cdots,a_\mathrm{k}^*\},$ for any $\xi \in \mathbb{Z}$. Also, there exists an absolute constant $C$ such that
$$|b(\xi)| \leq  C\langle \xi\rangle^{1-e}, \ \forall \xi \in \mathbb{Z}, \ and \ e > 0.$$

In order  to state our main results, we rewritten the wave equation \eqref{wave1} as new form, by introducing the new variables,
$$q=\mathrm{D}^{\frac{1}{2}}u+\mathrm{i}\mathrm{D}^{-\frac{1}{2}}\partial_tu, \quad \bar{q}=\mathrm{D}^{\frac{1}{2}}u-\mathrm{i}\mathrm{D}^{-\frac{1}{2}}\partial_tu,$$
where $$\mathrm{D}=\sqrt{-\Delta+\mathrm{m}}.$$
In the new variables , the equation \eqref{wave1} is transformed to
$$\mathrm{i}\partial_tq(t)=\mathrm{D}q(t)+\frac{1}{2}\mathrm{D}^{-\frac{1}{2}}\mathcal{W}(\omega t)\mathrm{D}^{-\frac{1}{2}}(q(t)+\bar{q}(t)).$$
Taking its complex conjugate, we can obtain the following matrix valued, Schr\"odinger-type system
\begin{equation}\label{2}
\mathrm{i}\partial_t\mathbf{q}(t)=\mathbf{H}(t)\mathbf{q}(t), \quad \mathbf{H}(t)= \mathbf{D}+ \epsilon \mathbf{K}(\omega t),
\end{equation}
\begin{equation}
\mathbf{D}=\left(
\begin{array}{cc}
\mathrm{D}& 0 \\
0 & -\mathrm{D}
\end{array}
\right ),\quad \mathbf{K}(\omega t)=\left(
\begin{array}{cc}
\mathcal{K}(\omega t)& \mathcal{K}(\omega t) \\
-\mathcal{K}(\omega t) & -\mathcal{K}(\omega t)
\end{array}
\right ),
\end{equation}
where $\mathcal{K}(\omega t)=\frac{1}{2}\mathrm{D}^{-\frac{1}{2}} \mathcal{W}(\omega t)\mathrm{D}^{-\frac{1}{2}} $ .

\begin{thm}\label{main1}
Assume the conditions $\mathrm{I}, \mathrm{II}$ of linear equation \eqref{2}. There exists $\bar{s}>0$, such that for any $s\geq \bar{s}$ there exists $\epsilon_{0}:=\epsilon(s,d)>0$ , $\gamma:=\gamma(s,d)>0$
and $\mathfrak{S}_s:=\mathfrak{S}(s,d)>0$ with $0\leq \mathfrak{S}_s\leq s$ such that if $\mathcal{W} \in \mathcal{OPS}^{1}_{s,2}$ satisfies
$$|\mathcal{W}|^{\gamma}_{1,s,2} \leq \epsilon,$$
then for any $\epsilon \in(0,\epsilon_{0})$ there exists a cantor like set $\mathcal{O}_{\epsilon} \in \mathcal{O}$ of asymptotically full Lebesgue
measure, i.e.
$$\lim_{\epsilon\mapsto 0}\mathrm{meas}(\mathcal{O} \setminus \mathcal{O}_{\epsilon}  )=0,$$
such that the following hold true. For any $\omega \subseteq\mathcal{O}_{\epsilon}$, there exists a liner bounded and invertible operator $\mathcal{T}(\omega t,\omega ) \in \mathcal{B}(\mathcal{H}^r_x\times \mathcal{H}^r_x)$ for any $0 \leq r \leq \mathfrak{S}_s $, such that the  change of coordinates $\mathbf{q}=\mathcal{T}(\omega t,\omega )\mathbf{v}$ conjugates
the equation \eqref{2} to the block diagonal time-independent system
\begin{equation}\label{3}
\mathrm{i}\partial_t\mathbf{v}(t)=\mathbf{H}^{\infty}_0(t)\mathbf{v}(t), \ \mathbf{H}^{\infty}(t)=\left(
\begin{array}{cc}
\mathcal{H}^{\infty}_0& 0 \\
0 & -\overline{\mathcal{H}^{\infty}_0}
\end{array}
\right ), \ \mathcal{H}_0^{\infty}=\mathrm{diag}\{\mathbf{h}^{\infty}_j| j\in \mathbb{N}\}
\end{equation}
$\mathbf{h}^{\infty}_0$ is a real number close to $\sqrt{\mathrm{m}}$, $\{\mathbf{h}^{\infty}_j\}_{j\neq 0}$ are $2\times 2$ self-adjoint matrices.
\end{thm}
\begin{rem}From Lemma \ref{estpd}, we known that  $\mathcal{K}(\omega t) \in \mathcal{OPS}^{-1}_{s-3,2}$  is a real and self-adjoint operator, and  satisfies
$$|\mathcal{K}|^{\gamma}_{-1,s-3,2} \lesssim_{s}\epsilon .$$
\end{rem}
\begin{cor}\label{main2}
For any $0\leq r \leq \mathfrak{R}_s$ and $\omega \in \mathcal{O}_{\epsilon}$, the solutions $\mathbf{q}(t, x):=(q(t,x),\bar{q}(t,x)) $ of equation \eqref{2} with initial
condition $\mathbf{q}(0, x):=(q(0,x),\bar{q}(0,x)) \in \mathcal{H}^r_x \times \mathcal{H}^r_x$ satisfies
$$c_r\|\mathbf{q}(0,x)\|_{\mathcal{H}^r_x\times\mathcal{H}^r_x} \leq \|\mathbf{q}(t,x)\|_{\mathcal{H}^r_x\times \mathcal{H}^r_x} \leq C_r\|\mathbf{q}(0,x)\|_{\mathcal{H}^r_x \times \mathcal{H}^r_x},$$
for some absolute constants $c_r,C_r >0$.
\end{cor}

\section{Linear operators}
\subsection{Matrix representation of linear operator}

Consider a linear operator $\mathcal{A}:L^2(\mathbb{T})\mapsto L^2(\mathbb{T})$, it action on a function $u \in L^2(\mathbb{T})$ as
$$\mathcal{A}[u]=\sum_{j,k\in\mathbb{Z}}\mathcal{A}^{k}_j\hat{u}_j e^{\mathrm{i}jx},$$
where $$\mathcal{A}^{k}_j=\frac{1}{2\pi}\int_{\mathbb{T}}\mathcal{A}[e^{\mathrm{i}jx}]e^{-\mathrm{i}kx}dx, \  \forall j,k \in \mathbb{Z}.$$
Also, we can identify the linear operator $\mathcal{A}$ as a infinite-dimensional block matrix
\begin{equation}\label{block}
\big( [\mathcal{A}]^{\beta}_{\alpha}\big)_{\alpha,\beta \in \mathbb{N}}
\end{equation}
where
$$[\mathcal{A}]^{\beta}_{\alpha}=\big(\mathcal{A}^{k}_{j}\big)_{|j|=\alpha, |k|=\beta}.$$
Note that the matrix $[\mathcal{A}]^{\beta}_{\alpha}$ is a linear operator from $\mathbb{E}_{\alpha}$ to $\mathbb{E}_{\beta}$, where
$$\mathbb{E}_{\alpha}:=\mathbf{Span}\{e^{\mathbf{i}jx}, j\in \mathbb{Z}, |j|=\alpha\}.$$
We also consider a smooth $\theta-$ dependent families of linear operator $\theta\mapsto \mathcal{A}(\theta),\mathbb{T}^d \mapsto \mathcal{B}(L^2(\mathbb{T}))$, which can be expanded  in Fourier series as
\begin{equation}
\mathcal{A}(\theta):=\sum_{\ell \in \mathbb{Z}^d}\hat{\mathcal{A}}(\ell)e^{\mathrm{i}\ell\cdot \theta}, \quad \hat{\mathcal{A}}(\ell):=\frac{1}{(2\pi)^d}\int_{\mathbb{T}^d}\mathcal{A}(\theta)e^{-\mathrm{i}\ell \cdot \theta}d\theta.
\end{equation}
It action on a function $u\in L^2(\mathbb{T}^d \times \mathbb{T})$ as
\begin{equation}\label{matre}
\mathcal{A}(\theta)[u(\theta,x)]=\sum_{\stackrel{j,k\in \mathbb{Z}}{\ell,\ell'\in \mathbb{Z}^d}}\hat{\mathcal{A}}(\ell-\ell')^k_j\hat{u}_j(\ell')e^{\mathrm{i}\ell\cdot\theta}e^{\mathrm{i}j\cdot x}.
\end{equation}
\begin{rem}\label{link}
From the matrix representation \eqref{matre}, we can regard the linear operator $\mathcal{A}$ as a pseudo-differential operator $\mathrm{Op}(a(\theta,x,\xi))$, for any $\xi=j \in \mathbb{Z} $, one has
\begin{equation}
a(\theta,x,j)=\sum_{k \in \mathbb{Z}}\mathcal{A}^k_j(\theta)e^{\mathrm{i}(k-j)x}.
\end{equation}
\end{rem}
\begin{defn}\label{matrixdefn}\textbf{(Matrix block decay norm)}
Let $\mathcal{A}$ be a $\theta-$dependent linear operator, $ \mathcal{A}(\theta):\mathbb{T}^d\mapsto\mathcal{B}(L^2(\mathbb{T}))$. Given $s\geq 0$, we say that $\mathcal{A} \in \mathcal{M}_s$, if and only if
\begin{equation}
|\mathcal{A}|_s:=\sup_{\alpha \in \mathbb{N}}\Big(\sum_{\stackrel{\ell \in \mathbb{Z}}{\beta \in \mathbb{N}}} \langle \ell, \beta-\alpha \rangle^{2s}\|[\hat{\mathcal{A}}(\ell)]^{\beta}_{\alpha}\|^2_0 \Big)^{\frac{1}{2}}< +\infty.
\end{equation}
$\|\cdot\|_0$ stands  for $L^2$ operator norm. If the operator $\mathcal{A}$ is Lipschitz depending on the parameters $\omega \in \mathcal{O}$. For any $\gamma \geq 0$. we claim that $\mathcal{A}(\omega) \in \mathrm{Lip}(\mathcal{O},\mathcal{M}_s)$, if and only if
\begin{equation}
|\mathcal{A}|^{\gamma,\mathcal{O}}_{s}:=\sup_{\omega \in \mathcal{O}}|\mathcal{A}(\omega)|_s+\gamma \sup_{\stackrel{\omega_1, \omega_2\in \mathcal{O}}{\omega_1 \neq \omega_2}}\frac{|\mathcal{A}(\omega_1)-\mathcal{A}(\omega_2)|_{s}}{|\omega_1-\omega_2|} < + \infty
\end{equation}
\end{defn}
\begin{rem}If the linear operator $\mathcal{A}$ is independent of $\theta$ and has the block matrix representation \eqref{block}, the matrix decay block norm becomes
\begin{equation}
|\mathcal{A}|_s:=\sup_{\alpha \in \mathbb{N}}\Big(\sum_{\beta \in \mathbb{N}} \langle  \beta-\alpha \rangle^{2s}\|[\mathcal{A}]^{\beta}_{\alpha}\|^2_0 \Big)^{\frac{1}{2}}.
\end{equation}
\end{rem}
\begin{lem}\label{est-mat}
$\mathbf{1:}$Let $s\geq s_0$, and $\mathcal{A}\in \mathcal{M}_s$,$u\in\mathcal{H}^s$. Then,
\begin{equation}
\|\mathcal{A}u\|_s\lesssim_{s}|\mathcal{A}|_s \|u\|_{s_0}+ |\mathcal{A}|_{s_0} \|u\|_s
\end{equation}

$\mathbf{2:}$ Let $s\geq s_0$, and $\mathcal{A},\mathcal{B}\in \mathcal{M}_s$. Then,
\begin{equation}
|\mathcal{AB}|_s \lesssim_s |\mathcal{A}|_s |\mathcal{B}|_{s_0}+ |\mathcal{A}|_{s_0} |\mathcal{B}|_{s_0}.
\end{equation}

$\mathbf{3:}$ Let $s\geq s_0$, and $\mathcal{A}\in \mathcal{M}_s$. There exists a constant $C(s)>0$, such that for any integer $n\geq 1$,
\begin{equation}
|\mathcal{A}^n|_s \leq C(s)^{n-1}(|\mathcal{A}|_{s_0})^{n-1}|\mathcal{A}|_s.
\end{equation}

$\mathbf{4:}$ Let $s\geq s_0$, and $\Phi:=exp(\mathcal{A})$ with $\mathcal{A} \in \mathcal{M}_s, \|\mathcal{A}\|_{s_0} \leq 1$. Then,
\begin{equation}
|\Phi^{\pm}-\mathrm{Id}|_s \lesssim_s|\mathcal{A}|_s.
\end{equation}

$\mathbf{5:}$ Items $\mathbf{1}-\mathbf{4}$ hold, replacing $|\cdot|_s$ by $|\cdot|^{\gamma}_s$ and $\|\cdot\|_s$ by $\|\cdot\|^{\gamma}_s$.

\end{lem}
\begin{proof}The proof is similar to Lemmata 2.7, 2.8 in \cite{Mon2017}.
\end{proof}
Given $N \in \mathbb{N}$, we define a smoothing operator $\Pi_{N}\mathcal{A}$ for any operator $\mathcal{A}$ with block -matrix representation \eqref{block},
$$
[\widehat{\Pi_N\mathcal{A}}(\ell)]^{\beta}_{\alpha}=\begin{cases}
[\hat{\mathcal{A}}(\ell)]^{\beta}_{\alpha}, \quad if \ |\ell| < N,\\
0, \quad \quad \quad \  otherwise.\\
\end{cases}$$

\begin{lem}\label{cut}
For any $s,\alpha>0$, the operator $\Pi^{\bot}_{N}:=\mathrm{Id}-\Pi_{N}$ stisfies
\begin{equation}
|\Pi^{\bot}_{N}\mathcal{A}|_s \leq N^{-\alpha}|\mathcal{A}|_{s+\alpha}, \quad |\Pi^{\bot}_{N}\mathcal{A}|^{\gamma}_s \leq N^{-\alpha}|\mathcal{A}|^{\gamma}_{s+\alpha}.
\end{equation}
\end{lem}
From the Definition \ref{defnpd},\ref{matrixdefn} and Remark \ref{link}, we can state a link between the pseudo-differential norms and matrix block decay norms.
\begin{lem}\label{link1}
Let $s\geq s_0$ and $\mathcal{A}\in \mathcal{OPS}^{0}_{s,0}$, one has
\begin{equation*}
|\mathcal{A}|^{\gamma}_{s} \lesssim_s |\mathcal{A}|^{\gamma}_{0,s,0}
\end{equation*}
\end{lem}
\begin{lem}\label{link2}(Lemma A.2 in \cite{F.G2020})
Let $\mathcal{A} \in \mathcal{B}(\mathcal{H}^r_{x})$ for $r\geq 0$ with $|\mathcal{A}|_{r+s_0} \leq +\infty$, then, one has
$$\sup_{\theta \in \mathbb{T}^d}\|\mathcal{A}(\theta)\|_{\mathcal{B}(\mathcal{H}^r_{x})} \lesssim_{r}|\mathcal{A}|_{r+s_0}.$$
\end{lem}
\begin{defn}\label{Ms}Let $\mathcal{A}$ be a $\theta-$dependent linear operator. Given $s\geq s_0 $ and $\rho > 0$, we say that $\mathcal{A} \in \mathcal{M}_{s,\rho}$ if
\begin{equation}
\langle D\rangle^{\rho}\mathcal{A},\quad  \mathcal{A}\langle D\rangle^{\rho} , \quad
\langle D\rangle^{\sigma}\mathcal{A} \langle D\rangle^{\sigma} \in \mathcal{M}_{s}, \quad \forall  \sigma\in\{0,\pm \rho\}.
\end{equation}
We can endow the $\mathcal{M}_{s,\rho}$ with the norm
\begin{equation}
\|\mathcal{A}\|_{s,\rho}:=|\langle D\rangle^{\rho}\mathcal{A}|_s+|\mathcal{A}\langle D\rangle^{\rho}|_s + \sum_{\sigma \in\{\pm\rho,0\}}|\langle D\rangle^{\sigma}\mathcal{A} \langle D\rangle^{-\sigma}|_s
\end{equation}
If the operator $\mathcal{A}$ is Lipschitz depending on the parameters $\omega \in \mathcal{O}$, we say that that $\mathcal{A}(\omega) \in \mathrm{Lip}(\mathcal{O},\mathcal{M}_{s,\rho})$. For any $\gamma >0$, we endow it with the norm
\begin{equation}
\|\mathcal{A}\|^{\gamma,\mathcal{O}}_{s,\rho}:=\sup_{\omega \in \mathcal{O}}\| \mathcal{A}(\omega)\|_{s,\rho}+\gamma \sup_{\stackrel{\omega_1, \omega_2\in \mathcal{O}}{\omega_1 \neq \omega_2}}\frac{\| \mathcal{A}(\omega_1)-\mathcal{A}(\omega_2)\|_{s,\rho}}{|\omega_1-\omega_2|}.
\end{equation}
\end{defn}
\begin{rem}For any $\sigma \in \mathbb{R}$, the operator $\langle D\rangle^{\sigma}$ is defined by $\langle D\rangle^{\sigma} e^{\mathrm{i}j\cdot x}=\langle j \rangle^{\sigma}e^{\mathrm{i}j\cdot x}$.
\end{rem}
\begin{defn}\label{Msp}Let $\mathcal{A}$ be a $\theta-$dependent linear operator. Given $s\geq s_0 $ and $\rho > 0$, we say that $\mathcal{A} \in \mathcal{L}_{s,\rho}$ if
\begin{equation}
\langle D\rangle^{\sigma}\mathcal{A} \langle D\rangle^{-\sigma} \in \mathcal{M}_{s}, \quad \forall  \sigma\in\{0,\pm \rho\}.
\end{equation}
We can endow the $\mathcal{L}_{s,\rho}$ with the norm
\begin{equation}
\|\mathcal{A}\|_{s}^{\rho}:=\sum_{\sigma \in\{\pm\rho,0\}}|\langle D\rangle^{\sigma}\mathcal{A} \langle D\rangle^{-\sigma}|_s
\end{equation}
If the operator $\mathcal{A}$ is Lipschitz depending on the parameters $\omega \in \mathcal{O}$, we say that that $\mathcal{A}(\omega) \in \mathrm{Lip}(\mathcal{O},\mathcal{L}_{s,\rho})$. For any $\gamma >0$, we endow it with the norm
\begin{equation}
\|\mathcal{A}\|_{s}^{\rho,\gamma,\mathcal{O}}:=\sup_{\omega \in \mathcal{O}}\| \mathcal{A}(\omega)\|_{s}^{\rho}+\gamma \sup_{\stackrel{\omega_1, \omega_2\in \mathcal{O}}{\omega_1 \neq \omega_2}}\frac{ \| \mathcal{A}(\omega_1)-\mathcal{A}(\omega_2)\|_{s}^{\rho} } {|\omega_1-\omega_2|}.
\end{equation}
\end{defn}

\begin{lem}\label{flow1}
Let $\rho >0$ and $s\geq s_0$. Assume that $\mathcal{A} \in \mathcal{M}_{s,\rho}$ and $\mathcal{B} \in \mathcal{L}_{s,\rho}$. Then, the following assertions hold true.\\
$\textbf{1}:$ For any $r \in [0,s-s_0]$ and $\theta \in \mathbb{T}^d$, the operator $e^{\mathrm{i}\mathcal{B}(\theta)} \in \mathcal{B}(\mathcal{H}^{r}_x )$ with the standard operator norm uniformly bounded in $\theta$.\\
$\textbf{2}:$ The commutator $\mathrm{i}[\mathcal{B},\mathcal{A}]:=\mathrm{i}(\mathcal{B}\mathcal{A}-\mathcal{A}\mathcal{B})$ belongs to $\mathcal{M}_{s,\rho}$ and satisfies
\begin{equation}\label{com1}
\|\mathrm{i}[\mathcal{B},\mathcal{A}]\|_{s,\rho} \lesssim_{s}\|\mathcal{A}\|_{s,\rho} \|\mathcal{B}\|^{\rho}_{s_0}+\|\mathcal{A}\|_{s_0,\rho} \|\mathcal{B}\|^{\rho}_{s}.
\end{equation}
$\textbf{3}:$ The operator $e^{\mathrm{i}\mathcal{B}}\mathcal{A}e^{-\mathrm{i}\mathcal{B}} $ belongs to $\mathcal{M}_{s,\rho}$ and with the quantitative bounds
\begin{equation}
 \| e^{\mathrm{i}\mathcal {B}}\mathcal{A}e^{-\mathrm{i}\mathcal{B}} \|_{s,\rho} \leq e^{2C_s \|\mathcal{B} \|^{\rho}_{s}} \|\mathcal{A}\|_{s,\rho},
\end{equation}

The analogous assertions hold true, if $\mathcal{A} \in \mathrm{Lip}(\mathcal{O},\mathcal{M}_{s,\rho})$ and $\mathcal{B} \in \mathrm{Lip}(\mathcal{O},\mathcal{M}_s)$.
\end{lem}

\begin{proof} The item \textbf{1} is a direct corollary of Lemmata \ref{est-mat}, \ref{link2} and Definition \ref{Ms}.

 For the commutator $\mathrm{i}[\mathcal{B},\mathcal{A}]:=\mathrm{i}(\mathcal{B}\mathcal{A}-\mathcal{A}\mathcal{B})$, the following inequalities hold (here $\sigma:=\pm \rho,0$)
 \begin{equation}
 \begin{split}
 |\langle D \rangle^{\rho}\mathcal{A}\mathcal{B}|_{s}&\lesssim_{s} |\langle D \rangle^{\rho}\mathcal{A}|_{s}|\mathcal{B}|_{s_0}+|\langle D \rangle^{\rho}\mathcal{A}|_{s_0}|\mathcal{B}|_{s},\\
 |\mathcal{A}\mathcal{B}\langle D \rangle^{\rho}|_{s}&\lesssim_{s} |\mathcal{A}\langle D \rangle^{\rho}|_{s}|\langle D \rangle^{-\rho}\mathcal{B}\langle D \rangle^{\rho}|_{s_0}+|\mathcal{A}\langle D \rangle^{\rho}|_{s_0}|\langle D \rangle^{-\rho}\mathcal{B}\langle D \rangle^{\rho}|_{s},\\
 |\langle D \rangle^{\rho}\mathcal{A}\mathcal{B}\langle D \rangle^{-\rho}|_{s} &\lesssim_{s} |\langle D \rangle^{\rho}\mathcal{A}\langle D \rangle^{-\rho}|_{s}|\langle D \rangle^{\rho}\mathcal{B}\langle D \rangle^{-\rho}|_{s_0}+|\langle D \rangle^{\rho}\mathcal{A}\langle D \rangle^{-\rho}|_{s_0}|\langle D \rangle^{\rho}\mathcal{B}\langle D \rangle^{-\rho}|_{s},
 \end{split}
 \end{equation}
 the same inequalities hold for $\mathcal{BA}$. Thus, one can get  the the quantitative bounds
 \begin{equation}
\|\mathrm{i}[\mathcal{A},\mathcal{B}]\|_{s,\rho} \lesssim_{s}\|\mathcal{A}\|_{s,\rho} \|\mathcal{B}\|^{\rho}_{s_0}+\|\mathcal{A}\|_{s_0,\rho} \|\mathcal{B}\|^{\rho}_{s}.
\end{equation}

For the operator $e^{\mathrm{i}\mathcal{B}}\mathcal{A}e^{-\mathrm{i}\mathcal{B}} $,  one has
\begin{equation}
e^{\mathrm{i}\mathcal{B}}\mathcal{A}e^{-\mathrm{i}\mathcal{B}} :=\mathcal{A}+\mathrm{i}[\mathcal{B},\mathcal{A}]+\frac{\mathrm{i}[\mathcal{B},\mathrm{i}[\mathcal{B},\mathcal{A}]]}{2!}+\cdots.
\end{equation}
From \eqref{com1}, one gets
\begin{equation}
\begin{split}
\|e^{\mathrm{i}\mathcal{B}}\mathcal{A}e^{-\mathrm{i}\mathcal{B}}\|_{s,\rho} &\leq \|\mathcal{A}\|_{s,\rho}+\|\mathrm{i}[\mathcal{B},\mathcal{A}]\|_{s,\rho}+\|\frac{\mathrm{i}[\mathcal{B},\mathrm{i}[\mathcal{B},\mathcal{A}]]}{2!}\|_{s,\rho}+\cdots\\
&\leq \|\mathcal{A}\|_{s,\rho}+2C_s\|\mathcal{A}\|_{s,\rho}\|\mathcal{B}\|_{s}^{\rho}+\frac{2^2C^2_s(\|\mathcal{B}\|^{\rho}_{s})^2\|\mathcal{A}\|_{s,\rho}}{2!}+\cdots\\
&\leq \|\mathcal{A}\|_{s,\rho}e^{2C_s \|\mathcal{B} \|^{\rho}_{s}}
\end{split}
\end{equation}

\end{proof}
\begin{lem}\label{opera1}
Let $s\geq s_0, \rho >0$ and $\mathcal{A}\in \mathcal{OPS}^{-\rho}_{s+\rho,0}$, one has
\begin{equation*}
\|\mathcal{A}\|^{\gamma}_{s,\rho} \lesssim_{s,\rho} |\mathcal{A}|^{\gamma}_{-\rho,s+\rho,0}.
\end{equation*}
Also, if  $\mathcal{A}\in \mathcal{OPS}^{0}_{s+\rho,0}$, one has
\begin{equation*}
\|\mathcal{A}\|^{\rho,\gamma}_{s} \lesssim_{s,\rho} |\mathcal{A}|^{\gamma}_{0,s+\rho,0}.
\end{equation*}
\end{lem}
\begin{proof}The proof is  a direct corollary of Lemma \ref{link1} and Definitions \ref{Ms}, \ref{Msp}.
\end{proof}

\subsection{The real and self-adjoint operators}
\begin{defn} $\mathbf{(1)}$: Given a $\theta$-dependent operator $\mathcal{A}(\theta):\mathbb{T}^d\mapsto \mathcal{B}(L^2(\mathbb{T}))$, we define its conjugate operator $\bar{\mathcal{A}}$ by $\bar{\mathcal{A}}u=\overline{\mathcal{A}\bar{u}}$. The conjugate operator $\bar{\mathcal{A}}$ has the matrix representation
$$\bar{\mathcal{A}}(\theta)^j_i=(\overline{\mathcal{A}(\theta)^{-j}_{-i}}),\quad \forall i,j \in \mathbb{Z},\ \theta \in \mathbb{T}^d.$$
$\mathbf{(2):}$ Given a $\theta$-dependent operator $\mathcal{A}(\theta):\mathbb{T}^d\mapsto \mathcal{B}(L^2(\mathbb{T}))$, we define its adjoint operator $\mathcal{A}^*(\theta)$ by
$$\int_{\mathbb{T}}\mathcal{A}(\theta)[u]\cdot \bar{v}dx=\int_{\mathbb{T}}u\cdot \overline{\mathcal{A}(\theta)^*[v]}dx.$$
The adjoint operator $\mathcal{A}^*$ has the matrix representation
$$\mathcal{A}^*(\theta)^j_i=\overline{\mathcal{A}(\theta)^i_j}, \quad \forall i,j \in \mathbb{Z},\ \theta \in \mathbb{T}^d. $$
\end{defn}

\begin{lem}$\mathbf{(1)}$: A linear operator $\mathcal{A}(\theta)$ is real, if it maps real functions to real and $\bar{\mathcal{A}}(\theta)=\mathcal{A}(\theta).$ Then, for any $i,j \in \mathbb{Z}, \alpha,\beta \in \mathbb{N}, \theta \in \mathbb{T}^d$, the following assertions hold true,
\begin{equation}
\mathcal{A}(\theta)_i^j=(\overline{\mathcal{A}(\theta)^{-j}_{-i}})\Leftrightarrow \hat{\mathcal{A}}(\ell)^j_i=\overline{\hat{\mathcal{A}}(-\ell)^{-j}_{-i}}\Leftrightarrow [\hat{\mathcal{A}}(\ell)]^{\beta}_{\alpha}= \overline{[\hat{\mathcal{A}}(-\ell)]^{\beta}_{\alpha}}.
\end{equation}

$\mathbf{(2)}$: A linear operator $\mathcal{A}(\theta)$ is self-adjoint, if $\mathcal{A}(\theta)=\mathcal{A}^*(\theta)$. Then, for any $i,j \in \mathbb{Z}, \alpha,\beta \in \mathbb{N}, \theta \in \mathbb{T}^d$, the following assertions hold true,
\begin{equation}
\mathcal{A}(\theta)_i^j=(\overline{\mathcal{A}(\theta)^{i}_{j}})\Leftrightarrow \hat{\mathcal{A}}(\ell)^j_i=\overline{\hat{\mathcal{A}}(-\ell)^{i}_{j}}\Leftrightarrow [\hat{\mathcal{A}}(\ell)]^{\beta}_{\alpha}= ([\hat{\mathcal{A}}(-\ell)]^{\alpha}_{\beta})^*.
\end{equation}

\end{lem}

\begin{lem}
Let $\mathcal{A}=\mathrm{Op}(a)\in \mathcal{OPS}^{m}_{s,\alpha}$. Then, $\mathcal{A}$ is real if and only if $a(\theta,x,\xi)= \overline{a(\theta,x,-\xi)}$.
\end{lem}

\subsection{$2 \times 2 $ operator matrix}

In this section, we will describe a special structure of operator matrix.
\begin{equation}\label{operator}
\mathbf{A}=\left(
\begin{array}{cc}
\mathcal{A}^{d}& \mathcal{A}^a \\
 -\overline{\mathcal{A}^a} & -\overline{\mathcal{A}^{d}}
\end{array}
\right)
\end{equation}
It needs emphasize that the diagonal operator $\mathcal{A}^d$ and anti-diagonal operator $\mathcal{A}^a$ have different symmetry properties and matrix decay properties. For the details, we show it in the following definition.

\begin{defn}\label{spaceM}
Given a $2\times 2$ operator matrix $\mathbf{A}$ of the form \eqref{operator}, and $\rho,o \in \mathbb{R},s\geq s_0$, we say that $\mathbf{A}$ belongs to $\mathcal{N}_{s}(\rho,o)$ if
\begin{equation}
[\mathcal{A}^d]^*=\mathcal{A}^d, \quad \quad [\mathcal{A}^a]^*=\overline{\mathcal{A}^a}
\end{equation}
and
\begin{equation}\label{dp}
\mathcal{A}^d  \in \mathcal{M}_{s,\rho} \quad  \mathcal{A}^a \in  \mathcal{M}_{s,o}
\end{equation}
We can endow the $\mathcal{N}_{s}(\rho,o)$ with the norm
\begin{equation}
\begin{split}
\interleave  \mathbf{A} \interleave_{s,\rho,o}:=&|\langle D\rangle^{\rho}\mathcal{A}^d|_s+|\mathcal{A}^d\langle D\rangle^{\rho}|_s + |\langle D\rangle^{o}\mathcal{A}^a|_s+|\mathcal{A}^a\langle D\rangle^{o}|_s \\
\\
&+\sum_{\substack{\sigma \in\{\pm\rho,0\} \\ \delta\in \{d,a\}}}|\langle D\rangle^{\sigma}\mathcal{A}^{\delta} \langle D\rangle^{-\sigma}|_s.
\end{split}
\end{equation}
If the operator $\mathcal{A}^d$ and $\mathcal{A}^o$ are  Lipschitz depending on the parameters $\omega \in \mathcal{O}$, we say that that $\mathbf{A}(\omega) \in \mathrm{Lip}(\mathcal{O},\mathcal{N}_s(\rho,o))$. For any $\gamma >0$, we endow it with the norm
\begin{equation}
\interleave \mathbf{A}\interleave^{\gamma,\mathcal{O}}_{s,\rho,o}:=\sup_{\omega \in \mathcal{O}}\interleave \mathcal{A}(\omega)\interleave_{s,\rho,o}+\gamma \sup_{\stackrel{\omega_1, \omega_2\in \mathcal{O}}{\omega_1 \neq \omega_2}}\frac{\interleave \mathbf{A}(\omega_1)-\mathbf{A}(\omega_2)\interleave_{s,\rho,o}}{|\omega_1-\omega_2|}.
\end{equation}
\end{defn}
\begin{rem}\

$\bullet$ The symmetry properties of the space $\mathcal{N}_s(\rho,o)$ is equivalent to ask that the operator matrix $\mathbf{A}$ is the Hamiltonian vector field of a real valued quadratic Hamilton. For the details, we refer to \cite{Mon2017, Mo2019}.\

$\bullet$  The decay properties \eqref{dp} of the diagonal operator $\mathcal{A}^d$  is essential to the measure estimate in section 4, we will show that the decay properties of the diagonal operator can be maintained in KAM iteration.
\end{rem}
\begin{lem}\label{com2}Let $\alpha >0$ and $s\geq s_0$. Assume that $\mathbf{A} \in \mathcal{N}_{s}(\rho,0)$ and $\mathbf{B} \in \mathcal{N}_{s}(\rho,\rho)$. Then $\mathrm{i}[\mathbf{B},\mathbf{A}]$ belongs to $ \mathcal{N}_{s}(\rho,\rho)$ and satisfies
\begin{equation}
\interleave \mathrm{i}[\mathbf{B},\mathbf{A}] \interleave_{s,\rho,\rho} \lesssim_s \interleave \mathbf{A}\interleave_{s,\rho,0 }\interleave \mathbf{B}\interleave_{s_0,\rho,\rho }+ \interleave \mathbf{A}\interleave_{s_0,\rho,0 }\interleave \mathbf{B}\interleave_{s,\rho,\rho }
\end{equation}
If $\mathbf{A} \in \mathrm{Lip}(\mathcal{O},\mathcal{N}_s(\rho,0))$ and $\mathbf{B} \in \mathrm{Lip}(\mathcal{O},\mathcal{N}_s(\rho,\rho))$. Then $\mathrm{i}[\mathbf{B},\mathbf{A}]$ belongs to $\mathrm{Lip}(\mathcal{O},\mathcal{N}_s(\rho,\rho))$ and satisfies
\begin{equation}
\interleave \mathrm{i}[\mathbf{B},\mathbf{A}] \interleave^{\gamma}_{s,\rho,\rho} \lesssim_s \interleave \mathbf{A}\interleave^{\gamma}_{s,\rho,0 }\interleave \mathbf{B}\interleave^{\gamma}_{s_0,\rho,\rho }+\interleave \mathbf{A}\interleave^{\gamma}_{s_0,\rho,0 }\interleave \mathbf{B}\interleave^{\gamma}_{s,\rho,\rho }.
\end{equation}
\end{lem}
\begin{proof} The proof is similar to Lemma 2.22 in \cite{L.M2019} and Lemma \ref{flow1}.
\end{proof}
\begin{lem}\label{MS2}
Let $\rho >0$ and $s\geq s_0$. Assume that $\mathbf{A} \in \mathcal{N}_{s}(\rho,0)$ and $\mathbf{B} \in \mathcal{N}_{s}(\rho,\rho)$. Then, the following assertions hold true.\\
$\mathbf{1}:$ For any $r \in [0,s]$ and $\theta \in \mathbb{T}^d$, the operator $e^{\mathrm{i}\mathcal{B}(\theta)} \in \mathcal{B}(\mathcal{H}^{r}_x\times\mathcal{H}^{r}_x )$ with the standard operator norm uniformly bounded in $\theta$.\\
$\mathbf{2}:$ The operator $e^{\mathrm{i}\mathbf{B}}\mathbf{A}e^{-\mathrm{i}\mathbf{B}} $ belongs to $\mathcal{N}_{s}(\rho,0)$ and $e^{\mathrm{i}\mathbf{B}}\mathbf{A}e^{-\mathrm{i}\mathbf{B}}-\mathbf{A} $ belongs to $\mathcal{N}_{s}(\rho,\rho)$ with the quantitative bounds
\begin{equation}
\interleave e^{\mathrm{i}\mathbf{B}}\mathbf{A}e^{-\mathrm{i}\mathbf{B}} \interleave_{s,\rho,0} \leq e^{2C_s \interleave\mathbf{B} \interleave_{s,\rho,\rho}} \interleave \mathbf{A}\interleave_{s,\rho,0},
\end{equation}
\begin{equation}
\interleave e^{\mathrm{i}\mathbf{B}}\mathbf{A}e^{-\mathrm{i}\mathbf{B}}-\mathbf{A} \interleave_{s,\rho,\rho} \leq 2C_s e^{2C_s \interleave\mathbf{B} \interleave_{s,\rho,\rho}} \interleave \mathbf{A}\interleave_{s,\rho,0}\interleave\mathbf{B} \interleave_{s,\rho,\rho}.
\end{equation}

The analogous assertions hold true, if $\mathcal{A} \in \mathrm{Lip}(\mathcal{O},\mathcal{N}_s(\rho,0))$ and $\mathcal{B} \in \mathrm{Lip}(\mathcal{O},\mathcal{N}_s(\rho,\rho))$.

\end{lem}
\begin{proof} The proof is similar to Lemma 2.23 in \cite{L.M2019} and Lemma \ref{flow1}.
\end{proof}
\section{Regularization procedure}

The goal  of this section is smoothing the diagonal part of the $2\times 2$ operator matrix $\mathbf{K}(\omega t)$ in equation \eqref{2}. In the following lemmas, we will conjugate  the diagonal parts of the perturbation  into  smoothing operators, while the anti-diagonal parts remain  bounded operators. The results of this section is essential to the KAM iteration in the following section. For any  $\gamma \in (0,1)$ and $\tau_0>d$, we introduce the set
\begin{equation}\label{non1}
\mathcal{O}_{\gamma}=\big\{\omega \in \mathcal{O}:|\omega \cdot \ell+j|\geq \frac{\gamma}{\langle \ell \rangle^{\tau_0}},  \quad \forall (\ell,j) \in \mathbb{Z}^{d+1}\backslash \{0\}\big\}.
\end{equation}
\begin{lem}\label{reg1} For any $\omega \in \mathcal{O}_{\gamma}$, and symbol $k \in \mathcal{S}^0_{s-3,2}$ with $k=k^*$, there exists a symbol $g\in \mathcal{S}^0_{s-\tau_0-5,2}$, with $g=g^*$, such that
\begin{equation}\label{sym1}
k-\langle k\rangle_{\theta,x}-\omega\cdot \partial_\theta g- \partial_x g \cdot \frac{\xi}{|\xi|}\chi(\xi) \in \mathcal{S}^{-1}_{s-\tau_0-7,2}.
\end{equation}
Furthermore, one has
\begin{equation}\label{re0}
|\mathrm{Op}(g)|^{\gamma}_{0,s-\tau_0-5,\alpha} \lesssim \frac{1}{\gamma}|\mathrm{Op}(k)|^{\gamma}_{0,s,2}
\end{equation}
and
$$|\mathrm{Op}(k-\langle k\rangle_{\theta,x}-\omega\cdot \partial_\theta g- \partial_x g \cdot \frac{\xi}{|\xi|}\chi(\xi))|^{\gamma}_{-1,s-\tau_0-7,2} \lesssim_s \frac{1}{\gamma} |\mathrm{Op}(k)|^{\gamma}_{0,s,2}.$$
\end{lem}
\begin{proof}Taking a cut-off function $\chi_1(\xi) \in \mathcal{C}^{\infty}(\mathbb{R},\mathbb{R})$, it is satisfying
\begin{equation}
\chi_1(\xi)=0, \quad  \forall|\xi| \leq \frac{3}{2}, \quad and   \quad \chi_1(\xi)=1, \quad \forall|\xi| \geq 2.
\end{equation}

For the symbol \eqref{sym1}, one has
\begin{equation}\label{sym2}
\begin{split}
&k-\langle k\rangle_{\theta,x}-\omega\cdot \partial_\theta g- \partial_x g \cdot \frac{\xi}{|\xi|}\chi(\xi)\\
&=\chi_1(k-\langle k\rangle_{\theta,x})-\omega\cdot \partial_\theta g- \partial_x g \cdot \frac{\xi}{|\xi|}\chi(\xi)+(1-\chi_1)(k-\langle k \rangle_{\theta,x}).
\end{split}
\end{equation}
From the Lemma \ref{adjsym1}, we known that $(1-\chi_1)(k-\langle k \rangle_{\theta,x}) \in \mathcal{S}^{-1}_{s,2}$ and
\begin{equation}\label{sym3}
|\mathrm{Op}((1-\chi_1)(k-\langle k \rangle_{\theta,x}))|^{\gamma}_{-1,s,2}\lesssim_{s}|\mathrm{Op}(k)|^{\gamma}_{0,s,2}.
\end{equation}
Our goal is determine a symbol $g$, such that
$$\chi_1(k-\langle k\rangle_{\theta,x})-\omega\cdot \partial_\theta g- \partial_x g \cdot \frac{\xi}{|\xi|}\chi(\xi) \in \mathcal{S}^{-1}_{s-\tau_0-5,2}. $$
Since we require that $g=g^*$, we look for a symbol of the form $g=\frac{q+q^*}{2}$. From \eqref{sym2}, one has
\begin{equation}\label{sym3}
\begin{split}
&\chi_1(k-\langle k\rangle_{\theta,x})-\omega\cdot \partial_\theta g- \partial_x g \cdot \frac{\xi}{|\xi|}\chi(\xi)\\
=&\chi_1(k-\langle k\rangle_{\theta,x})-\omega\cdot \partial_\theta q- \partial_x q \cdot \frac{\xi}{|\xi|}\chi(\xi)\\
&+(\omega \cdot \partial_{\theta}- \frac{\xi}{|\xi|}\chi(\xi)\cdot \partial_x)[\frac{q^*-q}{2}]
\end{split}
\end{equation}

Next, we look for a symbol $q$ that  satisfies
\begin{equation}
\chi_1(k-\langle k\rangle_{\theta,x})-\omega\cdot \partial_\theta q- \partial_x q \cdot \frac{\xi}{|\xi|}\chi(\xi)=0.
\end{equation}
For any $\omega \in \mathcal{O}_{\gamma}$, the symbol $q$ defined as
\begin{equation}
\begin{split}
q(\theta,x,\xi,\omega):=&(\omega \cdot \partial_{\theta}+\partial_x)^{-1}[k(\theta,x,\xi)-\langle k \rangle_{\theta,x}(\xi)]\chi^+_1(\xi)\\
&+(\omega \cdot \partial_{\theta}-\partial_x)^{-1}[k(\theta,x,\xi)-\langle k \rangle_{\theta,x}(\xi)]\chi^-_2(\xi),
\end{split}
\end{equation}
where $\chi^+_1(\xi):=\chi_1(\xi)\mathbb{I}_{\{\xi>0 \}},\chi^+_1(\xi):=\chi_1(\xi)\mathbb{I}_{\{\xi\leq 0 \}}$. $\mathbb{I}_{\{\xi>0 \}}(resp.\mathbb{I}_{\{\xi\leq 0 \}})$ is the characteristic function of the set $\{\xi \in \mathbb{R}: \xi >0\}(resp.\{\xi \in \mathbb{R}: \xi \leq 0\})$.

It is easy to verified that $\chi^+_1(\xi),\chi^-_1(\xi)$ are $\mathcal{C}^{\infty}$. From the non-resonance condition \eqref{non1} and Definition \ref{defnpd}, we known that
$(\omega \cdot \partial_{\theta}\pm \partial_x)^{-1}[k(\theta,x,\xi)-\langle k \rangle_{\theta,x}(\xi)] \in \mathcal{S}^0_{s-\tau_0-3,2}$,
\begin{equation}
|\mathrm{Op}\big((\omega \cdot \partial_{\theta}\pm \partial_x)^{-1}[\langle k \rangle_{\theta,x}(\xi)-k(\theta,x,\xi)]\big)|^{\gamma}_{0,s-\tau_0-3,2} \leq \frac{1}{\gamma}|\mathrm{Op}(k)|^{\gamma}_{0,s,2}.
\end{equation}
and
\begin{equation}
|\mathrm{Op}(q)|^{\gamma}_{0,s-\tau_0-3,2} \lesssim_{s} \frac{1}{\gamma}|\mathrm{Op}(k)|^{\gamma}_{0,s,2}.
\end{equation}

Note that  $k=k^*$, from Lemmata \ref{adjsym1}, \ref{adjsym2}, one gets
\begin{equation}
\begin{split}
|\mathrm{Op}(q)-\mathrm{Op}(q^*)|^{\gamma}_{-1,s-\tau_0-5,2} &\lesssim_{s}|\mathrm{Op}\big((\omega \cdot \partial_{\theta}\pm \partial_x)^{-1}[k(\theta,x,\xi)-\langle k \rangle_{\theta,x}(\xi)]\big)|^{\gamma}_{0,s-\tau_0-3,2}\\
&\lesssim_{s}\frac{1}{\gamma}|\mathrm{Op}(k)|^{\gamma}_{0,s,2}
\end{split}
\end{equation}
Hence, one gets
\begin{equation}\label{sym4}
|\mathrm{Op}\big((\omega \cdot \partial_{\theta}- \frac{\xi}{|\xi|}\chi(\xi)\partial_x)[\frac{q^*-q}{2}]\big)|^{\gamma}_{-1,s-\tau_0-6,2}\lesssim_{s}\frac{1}{\gamma}|\mathrm{Op}(k)|^{\gamma}_{0,s,2}.
\end{equation}

Combined \eqref{sym1},\eqref{sym3} and \eqref{sym4}, there exists a symbol $g=\frac{q+q^*}{2}$ such that
\begin{equation}
k-\langle k\rangle_{\theta,x}+\omega\cdot \partial_\theta g- \partial_x g \cdot \frac{\xi}{|\xi|}\chi(\xi)=(1-\chi_1)(k-\langle k \rangle_{\theta,x})+(\omega \cdot \partial_{\theta}- \frac{\xi}{|\xi|}\chi(\xi)\partial_x)[\frac{q^*-q}{2}],
\end{equation}
\begin{equation}
\begin{split}
|\mathrm{Op}(k-\langle k\rangle+\omega\cdot \partial_\theta g- \partial_x g \cdot \frac{\xi}{|\xi|}\chi(\xi))|^{\gamma}_{-1,s-\tau_0-6,2} &\leq |\mathrm{Op}((1-\chi_1)(k-\langle k \rangle))|^{\gamma}_{-1,s,2}\\
&\quad \ \ +|\mathrm{Op}\big((\omega \cdot \partial_{\theta}- \frac{\xi}{|\xi|}\chi(\xi)\partial_x)[\frac{q^*-q}{2}]\big)|^{\gamma}_{-1,s-\tau_0-6,2}\\
&\lesssim_{s}\frac{1}{\gamma}|\mathrm{Op}(k)|^{\gamma}_{0,s,2},
\end{split}
\end{equation}
and
\begin{equation}
|\mathrm{Op}(g)|^{\gamma}_{0,s-\tau_0-5,2}\lesssim_{s} \frac{1}{\gamma}|\mathrm{Op}(k)|^{\gamma}_{0,s-3,2}.
\end{equation}
\end{proof}

\begin{lem}\label{G1}Consider the linear equation \eqref{2} and assume conditions $\mathrm{I}, \mathrm{II}$. If the frequency vector $\omega \in \mathcal{O}_{\gamma}$, there exists a time dependent change of coordinates $$[q(t),\bar{q}(t)]^{T}=[e^{-\mathrm{i}\mathcal{G}(\omega t, \omega)}v(t), e^{\mathrm{i}\overline{\mathcal{G}}(\omega t, \omega)}\bar{v}(t)]^{T},$$where
$$\mathcal{G}(\omega t, \omega)=\mathrm{Op}(g(\theta,x,\xi,\omega))\in \mathcal{OPS}^0_{s-\tau_0-5,2},$$
that conjugates equation  \eqref{2} to
\begin{equation}\label{wave3}
\mathrm{i}\partial_t\mathbf{v}(t)=\widetilde{\mathbf{H}}(t)\mathbf{v}(t), \quad \widetilde{\mathbf{H}}(t)=\mathbf{H}_0+\mathbf{P}(\omega t,\omega),
\end{equation}
where
\begin{equation}
\mathbf{H}_0=\mathbf{D}+[\mathbf{K}], \quad \ \ [\mathbf{K}]=\left(
\begin{array}{cc}
\mathrm{Op}(\langle k \rangle_{\theta,x})& 0 \\
 0& -\mathrm{Op}(\langle k \rangle_{\theta,x})
\end{array}
\right )
\end{equation}
and
$$ \mathbf{P}(\omega t,\omega) \in \mathrm{Lip}(\mathcal{O}, \mathcal{N}_{s-\tau_0-11}(1,0)). $$ Also, there exists a constant $C_{s}$ depending on $s $ such that
$$
\interleave \mathbf{P}(\omega t,\omega) \interleave^{\gamma}_{s-\tau_0-11,1,0} \leq \frac{1}{\gamma} e^{2\frac{C_{s}}{\gamma}}|\mathcal{K}|^{\gamma}_{0,s-3,\alpha}.
$$
\end{lem}
\begin{proof}

First of all, we split the diagonal operator $\mathbf{D}$ in \eqref{2} into two parts, that is
\begin{equation}
\mathbf{D}=\mathbf{B}+ \mathbf{Z}, \quad \mathbf{B}=\left(
\begin{array}{cc}
\sqrt{-\Delta}& 0 \\
0 & -\sqrt{-\Delta}
\end{array}
\right ), \quad
\mathbf{Z}=\left(
\begin{array}{cc}
\mathcal{Z}& 0 \\
0 & -\mathcal{Z}
\end{array}
\right ).
\end{equation}
From Lemma 8.6 in \cite{S21}, we known that the operator $\mathcal{Z}$ is a real Fourier multiplies of order $-1$. Also, we divided  the perturbation  $\mathbf{K}(\omega t)$  into two parts, that is

\begin{equation}\label{K}
\mathbf{K}=\mathbf{K}_1+ \mathbf{K}_2, \quad \mathbf{K}_1=\left(
\begin{array}{cc}
\mathcal{K}& 0 \\
0 & -\mathcal{K}
\end{array}
\right ), \quad
\mathbf{K}_2=\left(
\begin{array}{cc}
0& \mathcal{K} \\
-\mathcal{K} & 0
\end{array}
\right ).
\end{equation}
Through the transformation $\mathbf{q}(t)=e^{-\mathrm{i}\mathbf{G}(\omega t,\omega)}\mathbf{v}(t)$, where
\begin{equation}\label{op1}
e^{-\mathrm{i}\mathbf{G}(\omega t,\omega)}=\left(
\begin{array}{cc}
e^{-\mathrm{i}\mathcal{G}(\omega t, \omega)}& 0 \\
0 & e^{\mathrm{i}\overline{\mathcal{G}}(\omega t, \omega)}
\end{array}
\right ),
\end{equation}
we can conjugated the equation \eqref{2} to
\begin{equation}
\mathrm{i}\partial_t\mathbf{v}(t)=\widetilde{\mathbf{H}}(t)\mathbf{v}(t),
\end{equation}
where
\begin{align}
\widetilde{\mathbf{H}}(t)=&e^{\mathrm{i}\mathbf{G}(\omega t,\omega)}\mathbf{H}(t)e^{-\mathrm{i}\mathbf{G}(\omega t,\omega)}-\int^1_0e^{\mathrm{i}s\mathbf{G}}\dot{\mathbf{G}}e^{-\mathrm{i}s\mathbf{G}}ds\\
\label{esti1}=&\mathbf{B}+\mathbf{Z}+\mathrm{i}[\mathbf{G},\mathbf{B}]+\mathbf{K}_1-\dot{\mathbf{G}}\\
\label{esti2}&-\mathrm{i}\int^1_0(1-s)e^{\mathrm{i}s\mathbf{G}}[\mathbf{G},\dot{\mathbf{G}}]e^{-\mathrm{i}s\mathbf{G}}ds\\
\label{esti3}&+\mathrm{i}\int^1_0e^{\mathrm{i}s\mathbf{G}}[\mathbf{G},\mathbf{K}_1]e^{-\mathrm{i}s\mathbf{G}}ds\\
\label{esti4}&+\mathrm{i}\int^1_0e^{\mathrm{i}s\mathbf{G}}[\mathbf{G},\mathbf{Z}]e^{-\mathrm{i}s\mathbf{G}}ds\\
\label{esti5}&+\mathrm{i}\int^1_0e^{\mathrm{i}s\mathbf{G}}[\mathbf{G},\mathrm{i}[\mathbf{G},\mathbf{B}]]e^{-\mathrm{i}s\mathbf{G}}ds\\
\label{esti6}&+e^{\mathrm{i}\mathbf{G}}\mathbf{K}_2e^{-\mathrm{i}\mathbf{G}}
\end{align}

The main goal is determined a $2\times 2$ operator matrix $\mathbf{G}$, where
\begin{equation}\label{G}
\mathbf{G}=\left(
\begin{array}{cc}
\mathcal{G}(\omega t, \omega)& 0 \\
0 & -\overline{\mathcal{G}}(\omega t, \omega)
\end{array}
\right ),\quad \mathcal{G}(\omega t, \omega)=\mathrm{Op}(g(\theta,x,\xi,\omega)),
\end{equation}
such that  $\mathrm{i}[\mathbf{G},\mathbf{B}]+\mathbf{K}_1-\dot{\mathbf{G}}$  becomes a smoothing operator matrix.  From \eqref{op1}, one gets
\begin{equation}
\mathrm{i}[\mathbf{G},\mathbf{B}]+\mathbf{K}_1-\dot{\mathbf{G}}=\left(
\begin{array}{cc}
\mathrm{i}[\mathcal{G},\sqrt{-\Delta}]-\omega\cdot\partial_{\theta} \mathcal{G}+\mathcal{K}& 0 \\
0 & \mathrm{i}[\bar{\mathcal{G}},\sqrt{-\Delta}]+\omega\cdot \partial_\theta \bar{\mathcal{G}}-\mathcal{K}
\end{array}
\right ),
\end{equation}
Take the operator $\mathrm{i}[\mathcal{G},\sqrt{-\Delta}]-\omega\cdot\partial_{\theta} \mathcal{G}+\mathcal{K}-\mathrm{Op}(\langle k\rangle_{\theta,x})$ as $\mathcal{P}^d_1$, from Lemma \ref{estpd} and Remark \ref{pse}, one gets
\begin{equation}
\begin{split}
\mathcal{P}^d_1+\mathrm{Op}(\langle k\rangle_{\theta,x})=&\mathrm{Op}\big(\{g,|\xi|\chi(\xi)\}\big)+\mathrm{i}\mathrm{Op}(r_{g,|\xi|\chi(\xi)})+\mathrm{Op}(-\omega\cdot \partial_{\theta}g)+\mathrm{Op}(k)\\
=&\mathrm{Op}(-\partial_xg \cdot \frac{|\xi|}{\xi}\chi(\xi)-\omega\cdot \partial_{\theta}g+k-\langle k \rangle_{\theta,x} )+\mathrm{Op}(-\partial_xg\cdot|\xi|\partial_{\xi}\chi(\xi))\\
&+\mathrm{i}\mathrm{Op}(r_{g,|\xi|\chi})+\mathrm{Op}(\langle k\rangle_{\theta,x}).
\end{split}
\end{equation}
From Lemma \ref{reg1}, there exists a symbol $g\in \mathcal{S}^0_{s-\tau_0-5,2}$ with $g=g^*$, such that
\begin{equation}\label{re1}
|\mathrm{Op}(k-\langle k\rangle-\omega\cdot \partial_\theta g- \partial_x g \cdot \frac{\xi}{|\xi|}\chi(\xi))|^{\gamma}_{-1,s-\tau_0-5,2} \lesssim_s \frac{1}{\gamma} |k|^{\gamma}_{0,s,2}.
\end{equation}
Since $\partial_{\xi}\chi(\xi)=0$ for $|\xi|\geq 1$, from Lemma \ref{adjsym1}, one gets
\begin{equation}\label{re2}
|\mathrm{Op}(-\partial_xg\cdot|\xi|\partial_{\xi}\chi(\xi))|^{\gamma}_{-1,s-\tau_0-6,2} \lesssim_{\alpha}|\mathrm{Op}(g)|^{\gamma}_{0,s-\tau_0-5,2}.
\end{equation}
From Lemma \ref{estpd} and Remark \ref{pse}, one gets
\begin{equation}\label{re3}
|\mathrm{Op}(r_{g,|\xi|\chi(\xi)})|^{\gamma}_{-1,s-\tau_0-10,0}\lesssim_{s}|\mathrm{Op}(g)|^{\gamma}_{0,s-\tau_0-5,2}.
\end{equation}
Combining the results of \eqref{re0}, \eqref{re1},\eqref{re2} and \eqref{re3}, from Lemma \ref{opera1}, we can get
\begin{equation}\label{p-1}
\|\mathcal{P}^d_1\|^{\gamma}_{s-\tau_0-11,1}\lesssim_s |\mathcal{P}^d_1|^{\gamma}_{-1,s-\tau_0-10,0}\lesssim_{s}\frac{1}{\gamma}|\mathrm{Op}(k)|^{\gamma}_{0,s-3,2}.
\end{equation}
Since  the operators $\mathcal{K,G},\sqrt{-\Delta}$ are self-adjoint,  the operator $\mathcal{P}^d_1$ is also self-adjoint. Note that the operators $\sqrt{-\Delta},\mathcal{K}$ are real, then
$$-\overline{\mathcal{P}^d_1}=\mathrm{i}[\bar{\mathcal{G}},\sqrt{-\Delta}]+\omega\cdot \partial_\theta \bar{\mathcal{G}}-\mathcal{K}+\mathrm{Op}(\langle k \rangle).$$
Finally, \eqref{esti1} can be rewritten as $\mathbf{H}_0+\mathbf{P}_1$, where
\begin{equation}
\mathbf{H}_0=\mathbf{B}+\mathbf{Z}+[\mathbf{K}],\quad  \mathbf{P}_1=\left(
\begin{array}{cc}
\mathcal{P}^d_1 & 0 \\
0 & -\overline{\mathcal{P}^d_1}
\end{array}
\right ).
\end{equation}
From \eqref{p-1}, one gets $\mathbf{P}_1 \in \mathcal{N}_{s-\tau_0-11}(1,0)$ and
\begin{equation}\label{P1}
 \interleave \mathbf{P}_1 \interleave^{\gamma}_{s-\tau_0-11,1,0} \lesssim_{s}\frac{1}{\gamma}|\mathrm{Op}(k)|^{\gamma}_{0,s-3,2}.
 \end{equation}
For the notational convenience, we rename \eqref{esti2}-\eqref{esti6} as $\mathbf{P}_2$ - $\mathbf{P}_6$.
The estimates of \eqref{esti2}-\eqref{esti5} are similar, so we take \eqref{esti3} as an example. From \eqref{K} and \eqref{G}, one sees  that $\mathbf{P}_3=\int^1_0\mathrm{i}e^{\mathrm{i}s\mathbf{G}}[\mathbf{G},\mathbf{K}_1]e^{-\mathrm{i}s\mathbf{G}}ds$, where
\begin{equation}
\begin{split}
\mathrm{i}e^{\mathrm{i}s\mathbf{G}}[\mathbf{G},\mathbf{K}_1]e^{-\mathrm{i}s\mathbf{G}}&=\left(
\begin{array}{cc}
e^{-\mathrm{i}\mathcal{G}(\omega t, \omega)}& 0 \\
0 & e^{\mathrm{i}\overline{\mathcal{G}}(\omega t, \omega)}
\end{array}\right )
\left(
\begin{array}{cc}
\mathrm{i}[\mathcal{G},\mathcal{K}]& 0 \\
0 & \mathrm{i}[\bar{\mathcal{G}},\mathcal{K}]
\end{array}
\right )
\left(
\begin{array}{cc}
e^{\mathrm{i}\mathcal{G}(\omega t, \omega)}& 0 \\
0 & e^{-\mathrm{i}\overline{\mathcal{G}}(\omega t, \omega)}
\end{array}
\right )\\
&=\left(
\begin{array}{cc}
\mathrm{i}e^{-\mathrm{i}\mathcal{G}(\omega t, \omega)}[\mathcal{G},\mathcal{K}]e^{\mathrm{i}\mathcal{G}(\omega t, \omega)}& 0 \\
0 &\mathrm{i} e^{\mathrm{i}\overline{\mathcal{G}}(\omega t, \omega)}[\bar{\mathcal{G}},\mathcal{K}]e^{-\mathrm{i}\overline{\mathcal{G}}(\omega t, \omega)}
\end{array}\right )
\end{split}
\end{equation}
From Lemmata \ref{estpd}, \ref{opera1} and Definitions \ref{Ms}, \ref{Msp}, one see that $\mathcal{G}\in \mathcal{L}_{s-\tau_0-6,1}$
\begin{equation}
\|\mathcal{G}\|^{1,\gamma}_{s-\tau_0-6}\lesssim_s |\mathcal{G}|^{\gamma}_{0,s-\tau_0-5,2} \lesssim_{s}\frac{1}{\gamma}|\mathcal{K}|^{\gamma}_{0,s-3,2}
\end{equation}
and
\begin{equation}
\|[\mathcal{G},\mathcal{K}]\|^{\gamma}_{s-\tau_0-9,1}\lesssim_s|[\mathcal{G},\mathcal{K}]|^{\gamma}_{-1,s-\tau_0-8,1}
\lesssim_{s}\frac{1}{\gamma}(|\mathcal{K}|^{\gamma}_{0,s-3,2})^2.
\end{equation}
From Lemma \ref{flow1}, one gets
\begin{equation}
\|\mathrm{i}e^{-\mathrm{i}\mathcal{G}}[\mathcal{G},\mathcal{K}]e^{\mathrm{i}\mathcal{G}}\|^{\gamma}_{s-\tau_0-9,1}\leq \frac{1}{\gamma} e^{2\frac{C_{s}}{\gamma}|\mathcal{K}|^{\gamma}_{0,s-3,2}}(|\mathcal{K}|^{\gamma}_{0,s-3,2})^2.
\end{equation}

Since $\mathrm{i}e^{-\mathrm{i}\mathcal{G}}[\mathcal{G},\mathcal{K}]e^{\mathrm{i}\mathcal{G}}$ is a self-adjoint operator, and $-\overline{\mathrm{i}e^{-\mathrm{i}\mathcal{G}}[\mathcal{G},\mathcal{K}]e^{\mathrm{i}\mathcal{G}}}=\mathrm{i} e^{\mathrm{i}\overline{\mathcal{G}}} [\bar{\mathcal{G}},\mathcal{K}] e^{-\mathrm{i}\overline{\mathcal{G}}}$, we can get $\mathbf{P}_3 \in \mathcal{N}_{s-\tau_0-6}(1,0)$ and
\begin{equation}\label{P3}
\interleave \mathbf{P}_3 \interleave^{\gamma}_{s-\tau_0-6,1,0} \leq \frac{1}{\gamma} e^{2\frac{C_{s}}{\gamma}|\mathcal{K}|^{\gamma}_{0,s-3,2}}(|\mathcal{K}|^{\gamma}_{0,s-3,2})^2.
\end{equation}

Repeat the estimate prcedure of $\mathbf{P}_3$ for $\mathbf{P}_2,\mathbf{P}_4,\mathbf{P}_5$, one can gets  $\mathbf{P}_2,\mathbf{P}_4,\mathbf{P}_5 \in \mathcal{N}_{s-\tau_0-11}(1,0)$ and
\begin{equation}\label{P2}
\interleave \mathbf{P}_2+\mathbf{P}_4+\mathbf{P}_5 \interleave^{\gamma}_{s-\tau_0-11,1,0} \lesssim \frac{1}{\gamma} e^{2\frac{C_{s}}{\gamma}|\mathcal{K}|^{\gamma}_{0,s-3,2}}|\mathcal{K}|^{\gamma}_{0,s-3,2}.
\end{equation}

For the operator matrix $\mathbf{P}_6$,  we can get
\begin{equation}
\begin{split}
\mathbf{P}_6&=e^{\mathrm{i}\mathbf{G}}\mathbf{K}_2e^{-\mathrm{i}\mathbf{G}}\\
&=\left(
\begin{array}{cc}
e^{-\mathrm{i}\mathcal{G}(\omega t, \omega)}& 0 \\
0 & e^{\mathrm{i}\overline{\mathcal{G}}(\omega t, \omega)}
\end{array}\right )
\left(
\begin{array}{cc}
0& \mathcal{K} \\
-\mathcal{K} & 0
\end{array}
\right )
\left(
\begin{array}{cc}
e^{\mathrm{i}\mathcal{G}(\omega t, \omega)}& 0 \\
0 & e^{-\mathrm{i}\overline{\mathcal{G}}(\omega t, \omega)}
\end{array}
\right )\\
&=\left(
\begin{array}{cc}
0 & e^{-\mathrm{i}\mathcal{G}}\mathcal{K} e^{-\mathrm{i}\overline{\mathcal{G}}} \\
-e^{\mathrm{i}\overline{\mathcal{G}}}\mathcal{K} e^{\mathrm{i}\mathcal{G}} &0
\end{array}\right )
\end{split}
\end{equation}
Take $e^{-\mathrm{i}\mathcal{G}}\mathcal{K} e^{-\mathrm{i}\overline{\mathcal{G}}}$ as $\mathcal{P}^a_6$, from Lemmata \ref{est-mat}, \ref{link1}, we can get
\begin{equation}
|\mathcal{P}^a_6|^{\gamma}_{s-\tau_0-6} \leq e^{2\frac{C_{s}}{\gamma}|\mathcal{K}|^{\gamma}_{0,s,\alpha}}|\mathcal{K}|^{\gamma}_{0,s,\alpha}.
\end{equation}
Since $[\mathcal{P}^a_6]^*= e^{\mathrm{i}\overline{\mathcal{G}}}\mathcal{K} e^{\mathrm{i}\mathcal{G}}=\overline{\mathcal{P}^a_6}$, we see that $\mathbf{P}_6 \in \mathcal{N}_{s-\tau-6}(1,0)$ and
\begin{equation}\label{P6}
\interleave \mathbf{P}_6 \interleave^{\gamma}_{s-\tau_0-6,1,0} \leq  e^{2\frac{C_{s}}{\gamma}|\mathcal{K}|^{\gamma}_{0,s-3,2}}|\mathcal{K}|^{\gamma}_{0,s-3,2}.
\end{equation}

Summing up $\mathbf{P}_1$ to $\mathbf{P}_6$, from \eqref{P1},\eqref{P3},\eqref{P2} and \eqref{P6}, one see that $\mathbf{P} \in \mathcal{N}_{s-\tau_0-11}(1,0)$ and
\begin{equation}\label{P}
\interleave \mathbf{P} \interleave^{\gamma}_{s-\tau_0-11,1,0} \lesssim \frac{1}{\gamma} e^{2\frac{C_{s}}{\gamma}}|\mathcal{K}|^{\gamma}_{0,s-3,2}.
\end{equation}
\end{proof}

\begin{rem}\label{regu}\

$\bullet$ For the non-resonance set $\mathcal{O}_{\gamma}$, it is easy to verified that
\begin{equation}
\mathrm{meas}(\mathcal{O}\backslash \mathcal{O}_{\gamma}) \leq C\gamma
\end{equation}
for an absolute constant.

$\bullet$ Take $\mathbf{V}^{\pm}(\omega t,\omega):=e^{\mp\mathrm{i}\mathbf{G}(\omega t,\omega)}$, from Lemma \ref{MS2}, the operator $\mathbf{V}^{\pm}(\omega t,\omega) $ belongs to $\mathcal{B}(\mathcal{H}^r_x \times \mathcal{H}^r_x )$ with standard operator norm for any $0 \leq r \leq s-\tau_0-11-s_0$.
\end{rem}

\section{KAM iteration}
After the preliminary transformation in the previous section, we can get the following equation
\begin{equation}\label{eqmain}
\mathrm{i}\partial_t\mathbf{v}(t)=\widetilde{\mathbf{H}}^0(t)\mathbf{v}(t), \quad \widetilde{\mathbf{H}}^0(t)=\mathbf{H}^0_0+\mathbf{P}^0(\omega t,\omega),
\end{equation}
where
\begin{equation}
\mathbf{P}^0=\mathbf{P},\quad \mathbf{H}^0_0=\mathbf{H}_0=\left(
\begin{array}{cc}
\mathcal{H}_0 & 0 \\
0 & -\mathcal{H}_0
\end{array}
\right ), \quad \mathcal{H}_0=\mathrm{diag}\{\mathbf{h}^0_j \ | \ j\in\mathbb{N}\},
\end{equation}
and $\mathbf{h}^0_j$ is linear operator from $\mathbb{E}_{j}$ to $\mathbb{E}_{j}$. To be more precise, for $j \in \mathbb{N}^+$ and $\mathfrak{a}\in \{1,-1\}$, one has
\begin{equation}
\mathbf{h}^j_0=\left(
\begin{array}{cc}
\lambda^0_{j,1} & 0 \\
0 & \lambda^0_{j,-1}
\end{array}
\right ), \quad \lambda^0_{j,\mathfrak{a}}=(j^2+\mathrm{m})^{\frac{1}{2}}+\langle k \rangle(\mathfrak{a}j).
\end{equation}
and $\mathbf{h}^0_0=\sqrt{\mathrm{m}}+\langle k \rangle_{\theta,x}(0)$.

\subsection{General step of KAM iteration}
In this section, we are going to perform an KAM iteration reducibility scheme for the linear equation \eqref{eqmain}. The main goal is to block-diagonalize the linear equation \eqref{eqmain}, and the key is to constantly square the size of the perturbation.

In the following, we show the outline of $k^{th}$ KAM iteration. For notional convenience, in the subsection, we drop the index $n$ and write $+$ instead of $n+1$.

Through a transformation $\mathbf{v}=e^{-\mathrm{i}\mathbf{U}^+}\mathbf{v}^{+}$, where
\begin{equation}\label{operator}
\mathbf{U}=\left(
\begin{array}{cc}
\mathcal{U}^{d}& \mathcal{U}^a \\
 -\overline{\mathcal{U}^a} & -\overline{\mathcal{U}^{d}}
\end{array}
\right),\quad [\mathcal{U}^{d}]^*=\mathcal{U}^{d}, \quad  [\mathcal{U}^{a}]^*=\overline{\mathcal{U}^{a}},
\end{equation}
the  equation $\mathrm{i}\partial_t\mathbf{v}(t)=\widetilde{\mathbf{H}}(t)\mathbf{v}(t)$ can be conjugated into
\begin{equation}
\mathrm{i}\partial_{t}\mathbf{v}^+=\widetilde{\mathbf{H}}^+(t)\mathbf{v}^+,
\end{equation}
where
\begin{align}
		\label{eq:reducibility_1}
		\widetilde{\mathbf{H}}^+(t) & = e^{\mathrm{i}\mathbf{U}^+(\omega t,\omega)}\widetilde{\mathbf{H}}(t)e^{-\mathrm{i}\mathbf{U}^+(\omega t,\omega)}-\int^1_0 e^{\mathbf{i}s\mathbf{U}^+(\omega t,\omega)}\dot{\mathbf{U}}^+e^{-\mathbf{i}s\mathbf{U}^+(\omega t,\omega)}ds, \\
		\label{eq:reducibility 2}
		& \widetilde{\mathbf{H}}^+= \mathbf{H}_0 +\mathbf{i}[\mathbf{U}^+,\mathbf{H}_0]+\mathbf{P} - \dot{\mathbf{U}}^+ +\mathbf{R},
	\end{align}
and
\begin{align}
		\label{eq:reducibility_3}
\mathbf{R} =& e^{\mathrm{i}\mathbf{U}^+(\omega t,\omega)}\mathbf{H}_0 e^{-\mathrm{i}\mathbf{U}^+(\omega t,\omega)}-(\mathbf{H}_0+\mathrm{i}[\mathbf{U}^+,\mathbf{H}_0])+(e^{\mathrm{i}\mathbf{U}^+(\omega t,\omega)}\mathbf{P}
e^{-\mathrm{i}\mathbf{U}^+(\omega t,\omega)}-\mathbf{P}) \\
&-(\int^1_0e^{\mathrm{i}s\mathbf{U}^+(\omega t,\omega)}\dot{\mathbf{U}}^+
e^{-\mathrm{i}s\mathbf{U}^+(\omega t,\omega)}ds-\dot{\mathbf{U}}^+  ).
\end{align}
We determine the operator matrix  $\mathbf{U}^+$ by solving the homological equation
\begin{equation}\label{hoeq}
\omega \cdot \partial_{\theta}\mathbf{U}^+=\mathrm{i}[\mathbf{U}^+,\mathbf{H}_0]+\Pi_{N}\mathbf{P}-\lfloor\mathbf{P}\rfloor,
\end{equation}
where
\begin{equation}\label{operator}
\lfloor\mathbf{P}\rfloor=\left(
\begin{array}{cc}
\lfloor \mathcal{P}^{d} \rfloor& 0 \\
 0 & -\overline{\lfloor \mathcal{P}^{d} \rfloor}
\end{array}
\right),\quad \lfloor \mathcal{P}^{d} \rfloor =\mathrm{diag}\{[\hat{\mathcal{P}^{d}}(0)]^j_j\ | \ j \in \mathbb{N}\}.
\end{equation}

The new Hamiltonian becomes $\mathbf{H}^{+}(t)=\mathbf{H}^{+}_0+\mathbf{P}^{+}$, where $\mathbf{P}^{+}=\mathbf{R}+ \Pi^{\bot}_{N}\mathbf{P}$ and
\begin{equation}\label{operator}
\mathbf{H}^{+}_0=\left(
\begin{array}{cc}
 \mathcal{H}^{+}_0 & 0 \\
 0 & -\overline{\mathcal{H}^{+}_0}
\end{array}
\right),\quad \mathcal{H}^{+}_0 =\mathrm{diag}\big \{\mathbf{h}^+_j \ \big | \ j \in \mathbb{N}\big \}=\mathrm{diag}\big \{\mathbf{h}_j+ [\hat{\mathcal{P}^{d}}(0)]^j_j\ \big | \ j \in \mathbb{N}\big \}.
\end{equation}

In order to get a nice solution of the homological equation \eqref{hoeq} and ensure the convergence of the KAM iteration, the second-order Melnikov conditions are required to be imposed.
Denoting $\lambda_{j,\mathfrak{a}},\mathfrak{a} \in \{1,-1\}$ as the eigenvalues of the block $\mathbf{h}_j$, we choose the
frequency vector $\omega$ from the following set
\begin{multline} \label{non-re}
\mathcal{O}^+_{\gamma}:=\Big\{\omega \in \mathcal{O}_{\gamma}:|\omega \cdot \ell+\lambda_{i,\mathfrak{a}}-\lambda_{j,\mathfrak{a}'}| \geq  \frac{\gamma \langle i- j\rangle}{N^{\tau}}, \quad \forall (\ell,i,j)\in \mathbb{Z}^d\times \mathbb{N} \times \mathbb{N},\ (\ell,i,j)\neq (0,j,j), \\ |\ell|\leq N, \
 and \ \{\mathfrak{a},\mathfrak{a}'\}\in \{-1,1\},\quad \quad
|\omega \cdot \ell+\lambda_{i,\mathfrak{a}}+\lambda_{j,\mathfrak{a}}| \geq  \frac{\gamma \langle i+ j\rangle}{N^{\tau}}, \\
\forall (\ell,i,j)\in \mathbb{Z}^d\times \mathbb{N} \times \mathbb{N},  \ |\ell|\leq N,\ \{\mathfrak{a},\mathfrak{a}'\}\in \{-1,1\}
\Big\}.
\end{multline}
\begin{lem}\label{h1}\textbf{(Homological Equation)} Assume that $ \mathbf{P}(\omega t,\omega) \in \mathrm{Lip}(\mathcal{O}_{\gamma}, \mathcal{N}_{s}(1,0))$, and
\begin{equation}\label{a1}
\max_{\omega \in \mathcal{O}}\frac{\|\Delta_{\omega}\mathbf{h}_j\|_{0}}{|\Delta \omega|} \leq C
\end{equation}
for an absolute constant. For any $\omega \in \mathcal{O}^+_{\gamma}$, there exists a solution $\mathbf{U}^+$ solve the equation \eqref{hoeq}, which belong to $\mathrm{Lip}(\mathcal{O}^+_{\gamma}, \mathcal{N}_{s}(1,1))$ with the quantitative bound
\begin{equation}
\interleave \mathbf{U}^+ \interleave^{\gamma}_{s,1,1} \lesssim \frac{N^{2\tau+1}}{\gamma}\interleave \mathbf{P}\interleave^{\gamma}_{s,1,0}.
\end{equation}
\end{lem}
\begin{proof}The homological equation \eqref{hoeq} is split in the two equation
\begin{equation}
-\mathrm{i}\omega\cdot\partial_{\theta}\mathcal{U}^d+[\mathcal{H}_0,\mathcal{U}^d]+\mathrm{i}\Pi_{N}\mathcal{P}^d=\mathrm{i}\lfloor\mathcal{P}^d\rfloor,
\end{equation}
\begin{equation}
-\mathrm{i}\omega\cdot\partial_{\theta}\mathcal{U}^a+\mathcal{H}_0\mathcal{U}^a+\mathcal{U}^a\overline{\mathcal{H}_0}+ \mathrm{i}\Pi_{N}\mathcal{P}^a=0.
\end{equation}

Considering the block matrix representation and expanding the Fourier series in time, for any $\ell \in \mathbb{Z}^d,|\ell|\leq N$ and $i,j\in \mathbb{N}$, we can get
\begin{equation}\label{hoeq1}
\omega \cdot \ell [\widehat{\mathcal{U}^d}(\ell)]^j_i+\mathbf{h}_i [\widehat{\mathcal{U}^d}(\ell)]^j_i- [\widehat{\mathcal{U}^d}(\ell)]^j_i\mathbf{h}_j=-\mathrm{i}[\widehat{\mathcal{P}^d}(\ell)]^j_i-[\lfloor\mathcal{P}^d\rfloor]^j_i,
\end{equation}
\begin{equation}\label{hoeq2}
\omega \cdot \ell [\widehat{\mathcal{U}^a}(\ell)]^j_i+\mathbf{h}_i [\widehat{\mathcal{U}^a}(\ell)]^j_i+ [\widehat{\mathcal{U}^a}(\ell)]^j_i\overline{\mathbf{h}}_j=-\mathrm{i}[\widehat{\mathcal{P}^d}(\ell)]^j_i.
\end{equation}
Takeing $\mathbf{I}$ as unite matrix, one see that, for any $i,j \in \mathbb{N}$, $$\mathrm{spec}(\omega\cdot\ell \mathbf{I}+\mathbf{h}_i)=\big\{\omega\cdot\ell+\lambda_{i,\mathfrak{a}} \ | \ \mathfrak{a}\in \{1,-1\} \big\}.$$
Since the block $\mathbf{h}_j$ is self-adjoint, one sees
$$\mathrm{spec}(\mathbf{h}_j)=\mathrm{spec}(\overline{\mathbf{h}}_j)=\big\{\omega\cdot\ell+\lambda_{j,\mathfrak{a}} \ | \ \mathfrak{a}\in \{1,-1\} \big\}.$$
From Lemma \ref{pt2}, one get immediately that
\begin{equation}\label{Ua}
\|[\widehat{\mathcal{U}^d}(\ell)]^j_i\|_0 \lesssim \frac{N^\tau{}}{\gamma}\|[\widehat{\mathcal{P}^d}(\ell)]^j_i\|_0.
\end{equation}
Let $\omega_1,\omega_2 \in \mathcal{O}^{+}_{\gamma}$, for any function $f=f(\omega)$, we write $\Delta_{\omega}f=f(\omega_1)-f(\omega_2)$. For the equation \eqref{hoeq1}, we can get that
\begin{equation}
\begin{split}
\omega \cdot \ell \Delta_{\omega}[\widehat{\mathcal{U}^d}(\ell)]^j_i+\mathbf{h}_i \Delta_{\omega}[\widehat{\mathcal{U}^d}(\ell)]^j_i- \Delta_{\omega}[\widehat{\mathcal{U}^d}(\ell)]^j_i\mathbf{h}_j&=-(\Delta \omega \cdot \ell) [\widehat{\mathcal{U}^d}(\ell)]^j_i(\omega_1)-\Delta_{\omega}\mathbf{h}_i [\widehat{\mathcal{U}^d}(\ell)]^j_i(\omega_1)\\
&+ [\widehat{\mathcal{U}^d}(\ell)]^j_i(\omega_1) \Delta_{\omega}\mathbf{h}_j-\mathrm{i}\Delta_{\omega}[\widehat{\mathcal{P}^d}(\ell)]^j_i-\mathrm{i}\Delta_{\omega}[\lfloor\mathcal{P}^d\rfloor]^j_i
\end{split}
\end{equation}
From Lemma \ref{pt2} again, we can get
\begin{equation}\label{Ua2}
\|\Delta_{\omega}[\widehat{\mathcal{U}^d}(\ell)]^j_i\|_0 \lesssim \frac{N^{\tau}}{\gamma}\|\Delta_{\omega}[\widehat{\mathcal{P}^d}(\ell)]^j_i\|_0+\frac{N^{2\tau+1}}{\gamma^2}\|[\widehat{\mathcal{P}^d}(\ell)]^j_i\|_0|\Delta \omega|
\end{equation}
Thus, \eqref{Ua},\eqref{Ua2} imply that
\begin{equation}\label{Ua5}
|\mathcal{U}^d|^{\gamma}_{s} \lesssim \frac{N^{\tau}}{\gamma} |\mathcal{P}^d|^{\gamma}_{s}.
\end{equation}
Also, the norms of $|\langle D \rangle\mathcal{U}^d|^{\gamma}_s$, $|\mathcal{U}^d\langle D \rangle|^{\gamma}_s$, $|\langle D \rangle^{\sigma}\mathcal{U}^d\langle D \rangle^{-\sigma}|^{\gamma}_s,(\sigma=\pm1)$ can be bounded by the same norms of $\mathcal{P}^d$.

The bounds control of $\mathcal{U}^a$ is more delicate. From Lemma \ref{pt2}, we can get that
\begin{equation}\label{Ua1}
\|[\widehat{\mathcal{U}^a}(\ell)]^j_i\|_0 \lesssim \frac{N^\tau}{\gamma\langle i+j\rangle}\|[\widehat{\mathcal{P}^a}(\ell)]^j_i\|_0.
\end{equation}
We also need control the bounds of $[\widehat{\langle D \rangle \mathcal{U}^a}(\ell)]^j_i,$ $[\widehat{ \mathcal{U}^a\langle D \rangle}(\ell)]^j_i $ and $[\Widehat{\langle D \rangle^{\sigma}\mathcal{U}^a\langle D \rangle^{-\sigma}}(\ell)]^j_i$, $(\sigma=\pm1)$. Considering the term $[\widehat{\langle D \rangle \mathcal{U}^a}(\ell)]^j_i$, and applying Lemma \ref{pt2} again, one gets
\begin{equation}\label{Ua3}
\begin{split}
\|[\widehat{\langle D \rangle \mathcal{U}^a}(\ell)]^j_i\|_0&=\|\langle i \rangle[\widehat{\mathcal{U}^a}(\ell)]^j_i\|_0 \lesssim\frac{N^{\tau}}{\gamma}\frac{\langle i \rangle}{\langle i+j \rangle}\|[\widehat{\mathcal{P}^a}(\ell)]^j_i\|_0\\
&\lesssim \frac{N^\tau}{\gamma}\|[\widehat{\mathcal{P}^a}(\ell)]^j_i\|_0.
\end{split}
\end{equation}
The similar bounds  hold for $[\widehat{ \mathcal{U}^a\langle D \rangle}(\ell)]^j_i $ and $[\Widehat{ \langle D \rangle^{\sigma}\mathcal{U}^a\langle D \rangle^{-\sigma}}(\ell)]^j_i$. Applying the difference operator $\Delta_{\omega }$ to equation \eqref{hoeq2}, one gets
\begin{equation}
\begin{split}
\omega \cdot \ell \Delta_{\omega}[\widehat{\mathcal{U}^a}(\ell)]^j_i+\mathbf{h}_i \Delta_{\omega}[\widehat{\mathcal{U}^a}(\ell)]^j_i+ \Delta_{\omega}[\widehat{\mathcal{U}^a}(\ell)]^j_i\overline{\mathbf{h}}_j&=-(\Delta \omega \cdot \ell) [\widehat{\mathcal{U}^a}(\ell)]^j_i(\omega_1)-\Delta_{\omega}\mathbf{h}_i [\widehat{\mathcal{U}^a}(\ell)]^j_i(\omega_1)\\
&+ [\widehat{\mathcal{U}^a}(\ell)]^j_i(\omega_1) \Delta_{\omega}\overline{\mathbf{h}}_j-\mathrm{i}\Delta_{\omega}[\widehat{\mathcal{P}^a}(\ell)]^j_i.
\end{split}
\end{equation}
Applying Lemma \ref{pt2} again, we can get
\begin{equation}\label{Ua4}
\| \Delta_{\omega}[\widehat{\mathcal{U}^a}(\ell)]^j_i\|_0  \lesssim \frac{N^{\tau}}{\gamma \langle i+j \rangle}\|\Delta_{\omega}[\widehat{\mathcal{P}^a}(\ell)]^j_i\|_0+\frac{N^{2\tau+1}}{\gamma^2\langle i+j \rangle^2}\|[\widehat{\mathcal{P}^d}(\ell)]^j_i\|_0|\Delta \omega|.
\end{equation}
The similar bounds  hold for $\Delta_{\omega}[\widehat{\langle D \rangle \mathcal{U}^a}(\ell)]^j_i,$ $\Delta_{\omega}[\widehat{ \mathcal{U}^a\langle D \rangle}(\ell)]^j_i $ and $\Delta_{\omega}[\Widehat{\langle D \rangle^{\sigma}\mathcal{U}^a\langle D \rangle^{-\sigma}}(\ell)]^j_i$, $(\sigma=\pm1)$.

Finally, \eqref{Ua5}, \eqref{Ua1},\eqref{Ua3},\eqref{Ua4} and Definition \ref{spaceM} imply that
\begin{equation}
\interleave \mathbf{U}^+ \interleave^{\gamma}_{s,1,1} \lesssim \frac{N^{2\tau+1}}{\gamma}\interleave \mathbf{P}\interleave^{\gamma}_{s,1,0}.
\end{equation}

\end{proof}

\begin{lem}\label{p1}\textbf{(The New Perturbation)}Fix $s\geq s_0$ and $b>0$. Let  $ \mathbf{P}(\omega t,\omega) \in \mathrm{Lip}(\mathcal{O}_{\gamma}, \mathcal{N}_{s+b}(1,0))$. Assume \eqref{a1} and,  for some fixed constant $C_{s}$,
\begin{equation}\label{as1}
C_{s}\frac{N^{2\tau+1}}{\gamma}\interleave\mathbf{P} \interleave^{\gamma}_{s_0,1,0} \leq \frac{1}{2}.
\end{equation}
Then, $\mathbf{P}^+=\Pi_{N}\mathbf{P}+\mathbf{R}$ is defined on $\mathcal{O}_{\gamma}^+$ and satisfies the quantitative bounds
\begin{equation}\label{PES1}
\interleave\mathbf{P}^+ \interleave^{\gamma}_{s,1,0} \lesssim_{s} N^{-b}\interleave\mathbf{P} \interleave^{\gamma}_{s+b,1,0} +\frac{N^{2\tau+1}}{\gamma}\interleave\mathbf{P} \interleave^{\gamma}_{s ,1,0}\interleave\mathbf{P} \interleave^{\gamma}_{s_0 ,1,0}
\end{equation}
\begin{equation}\label{PES2}
\interleave\mathbf{P}^+ \interleave^{\gamma}_{s+b,1,0} \leq C(s+b)\interleave \mathbf{P} \interleave^{\gamma}_{s+b,1,0}.
\end{equation}
\end{lem}
\begin{proof}From homological equation \eqref{hoeq}, the operator $\mathbf{R}$ defined in \eqref{eq:reducibility_3} can be rewritten as
\begin{equation}
\mathbf{R}=\mathrm{i}\int^1_0(1-s)e^{\mathrm{i}s\mathbf{U}^+}[\mathbf{U}^+,\lfloor\mathbf{P}\rfloor-\Pi_{N}\mathbf{P}]e^{-\mathrm{i}s\mathbf{U}^+}ds
+\mathrm{i}\int^1_0e^{\mathrm{i}s\mathbf{U}^+}[\mathbf{U}^+,\mathbf{P}]e^{-\mathrm{i}s\mathbf{U}^+}ds
\end{equation}
From Lemmata \ref{com2}, \ref{MS2}, \ref{h1}, one gets
\begin{equation}
\interleave \mathrm{i}[\mathbf{U}^+,\mathbf{P}]\interleave^{\gamma}_{s,1,1} \lesssim_{s}\frac{N^{2\tau+1}}{\gamma}\interleave \mathbf{P}\interleave^{\gamma}_{s,1,0}\interleave \mathbf{P}\interleave^{\gamma}_{s_0,1,0}
\end{equation}
and
\begin{equation}\label{PES}
\begin{split}
\interleave e^{\mathrm{i}\mathbf{U}^+}[\mathbf{U}^+,\mathbf{P}]e^{-\mathrm{i}\mathbf{U}^+}\interleave^{\gamma}_{s,1,1}& \leq \interleave \mathrm{i}[\mathbf{U}^+,\mathbf{P}]\interleave^{\gamma}_{s,1,1} + \interleave  \mathrm{i}[\mathbf{U}^+, \mathrm{i}[\mathbf{U}^+,\mathbf{P}]]\interleave^{\gamma}_{s,1,1}+\cdots\\
&\leq \frac{C_sN^{2\tau+1}}{\gamma}\interleave \mathbf{P}\interleave^{\gamma}_{s,1,0}\interleave \mathbf{P}\interleave^{\gamma}_{s_0,1,0}+ (\frac{2 C_sN^{2\tau+1}}{\gamma})^2\interleave \mathbf{P}\interleave^{\gamma}_{s,1,0}(\interleave \mathbf{P}\interleave^{\gamma}_{s_0,1,0})^2+\cdots \\
& \leq  e^{2C_s\frac{N^{2\tau+1}}{\gamma}\interleave \mathbf{P}\interleave^{\gamma}_{s_0,1,0}}\frac{N^{2\tau+1}}{\gamma}\interleave \mathbf{P}\interleave^{\gamma}_{s,1,0}\interleave \mathbf{P}\interleave^{\gamma}_{s_0,1,0}
\end{split}
\end{equation}
The same bounds hold for $\interleave \mathrm{i}[\mathbf{U}^+,\lfloor\mathbf{P}\rfloor-\Pi_{N}\mathbf{P}]\interleave^{\gamma}_{s,1,1}$ and $\interleave e^{\mathrm{i}\mathbf{U}^+}[\mathbf{U}^+,\lfloor\mathbf{P}\rfloor-\Pi_{N}\mathbf{P}]e^{-\mathrm{i}\mathbf{U}^+}\interleave^{\gamma}_{s,1,1}$.
From Lemma \ref{cut} and Definition \ref{spaceM}, one gets

\begin{equation}\label{PES3}
 \interleave \Pi_{N}\mathbf{P}\interleave^{\gamma}_{s,1,0} \leq N^{-b}\interleave \mathbf{P}\interleave^{\gamma}_{s,1,0}.
\end{equation}

Finally, \eqref{PES} and \eqref{PES3} imply \eqref{PES1}. By \eqref{as1}, the estimate of \eqref{PES2} can be obtained in the same way.

\end{proof}

\subsection{Iterative procedure}
The KAM iteration is start with the linear equation \eqref{eqmain}. The iteration objects are construct in Lemmata \ref{h1}, \ref{p1} by setting for $n\geq 0$
\begin{equation}
\begin{split}
\mathrm{i}\partial_{t}\mathbf{v}=\widetilde{\mathbf{H}}^{n}(t)\mathbf{v}&,\ \ \widetilde{\mathbf{H}}^{n}(t):=\mathbf{H}^{n}_0(\omega)+\mathbf{P}^{n}(\omega t,\omega), \  \ \mathbf{U}:=\mathbf{U}^{n}(\omega t ,\omega ),\\
&\mathbf{R}:=\mathbf{R}^{n}(\omega t ,\omega ), \ \  \mathbf{h}_j:= \mathbf{h}^{n}_j,\ \ \lambda_{j,\mathfrak{a}}:=\lambda_{j,\mathfrak{a}}^{n}.
\end{split}
\end{equation}

 We define
 \begin{equation}
 N_{-1}:=1, \quad N_{n}:=N^{(\frac{3}{2})^n}_0, \quad \forall n \geq 0,
 \end{equation}
 and for $\tau_0 >0, \tau>0$, we define the constants
 \begin{equation}\label{as1}
 \mathrm{s}:=s-\tau_0-11, \quad \mathrm{a}:=6\tau+4£¬\quad b= \mathrm{a}+1, \quad \mathrm{s}-b>s_0.
 \end{equation}
 Also, we assume that
 \begin{equation}
 \varepsilon:=\max \{\epsilon,  \interleave \mathbf{P} \interleave^{\gamma}_{s,1,0} \}.
 \end{equation}

 \begin{thm}\label{kam}\textbf{(KAM Reducibility)} Let $\gamma \in (0,1)$ and  $\mathrm{s}$ satisfies \eqref{as1}. There exists $N_0:=N_0(s,\tau,d) \in \mathbb{N}$ large enough and $\delta_0:=\delta_0(s,\tau,d) \in (0,1)$ such that if
 \begin{equation}
 \varepsilon \gamma^{-1} \leq \delta_0
 \end{equation}
 then:

 $\mathbf{I}:$ Seeing $\mathcal{O}^0_{\gamma}$ as  $\mathcal{O}_{\gamma}$ in \eqref{non1}, we can recursively defined for $n\geq 0$ and $\mathfrak{a},\mathfrak{a}' \in \{1,-1\}$,
 \begin{multline} \label{non-re}
\mathcal{O}^{n+1}_{\gamma}:=\Big\{\omega \in \mathcal{O}^n_{\gamma}:|\omega \cdot \ell+\lambda^n_{i,\mathfrak{a}}-\lambda^n_{j,\mathfrak{a}'}| \geq  \frac{\gamma \langle i- j\rangle}{N_n^{\tau}}, \quad \forall (\ell,i,j) \in \mathbb{Z}^d\times \mathbb{N} \times \mathbb{N}, \\ (\ell,i,j)\neq (0,j,j)
 ,\ \  |\ell|\leq N_n, \ and \
|\omega \cdot \ell+\lambda^n_{i,\mathfrak{a}}+\lambda^n_{j,\mathfrak{a}'}| \geq  \frac{\gamma \langle i+ j\rangle}{N^{\tau}}\\
 \forall (\ell,i,j) \in \mathbb{Z}^d\times \mathbb{N} \times \mathbb{N},\ (\ell,i,j)\neq (0,j,j), \ |\ell|\leq N_{n}
\Big\}.
\end{multline}

$\mathbf{II}:$ There exists a operator matrix  $\mathbf{U}^{n}(\omega t,\omega) \in \mathrm{Lip}(\mathcal{O}^{n}_{\gamma}, \mathcal{N}_{\mathrm{s}}(1,1))$ and satisfies
\begin{equation}
\interleave \mathbf{U}^{n} \interleave^{\gamma}_{s-b,1,1} \leq \interleave \mathbf{P}^{0} \interleave^{\gamma}_{s,1,0}\gamma^{-1}N_{n-1}^{2\tau+1}N_{n-2}^{-\mathrm{a}}.
\end{equation}
The change of coordinate $e^{-\mathrm{i}\mathbf{U}^{n}}$ conjugate $\widetilde{\mathbf{H}}^{n-1}$ to $\widetilde{\mathbf{H}}^{n}:=\mathbf{H}^{n}_0(\omega)+\mathbf{P}^{n}(\omega t,\omega)$ such that

$\mathbf{III}:$ The operator $\mathbf{H}^{n}_0$ is block diagonal, self-adjoint and time independent, where
\begin{equation}
\mathbf{H}^n_0=\left(
\begin{array}{cc}
\mathcal{H}^n_0 & 0 \\
0 & -\overline{\mathcal{H}^n_0}
\end{array}
\right ), \quad \mathcal{H}^n_0=\mathrm{diag}\{\mathbf{h}^n_j \ | \ j\in\mathbb{N}\},
\end{equation}
and the block $\mathbf{h}^n_j$ is defined over $\mathcal{O}_{\gamma}$, satisfies
\begin{equation}\label{h21}
\|\mathbf{h}^n_j-\mathbf{h}^{n-1}_j\|^{\gamma}_{0} \leq N^{-\mathrm{a}}_{n-2}\varepsilon j^{-1}.
\end{equation}

$\mathbf{IV}:$ The new perturbation $\mathbf{P}^{n}$ belongs to $\mathrm{Lip}(\mathrm{O}^{n}_{\gamma},\mathcal{N}_{\mathrm{s}}(1,0))$ and fulfils
\begin{equation}
\interleave \mathbf{P}^{n} \interleave^{\gamma}_{\mathrm{s},1,0} \leq \interleave \mathbf{P}^{0} \interleave^{\gamma}_{\mathrm{s},1,0}N_{n-1}.
\end{equation}
\begin{equation}
\interleave \mathbf{P}^{n} \interleave^{\gamma}_{\mathrm{s}-b,1,0} \leq \interleave \mathbf{P}^{0} \interleave^{\gamma}_{s,1,0}N^{-\mathrm{a}}_{n-1}.
\end{equation}
 \end{thm}
 \begin{proof}We prove these assertions by inductive. It's easy to verified that  the properties in items $\mathbf{I}-\mathbf{IV}$ hold true for $n=0$. Let us suppose that the statements hold true for a fixed $n\in \mathbb{N}$ and define the set $\mathcal{O}^{n+1}_{\gamma}$ in item $\mathbf{I}$. We prove these assertions also hold true for $n+1$.

 In order to apply  Lemmata \ref{h1}, \ref{p1}, we need check the assumptions in the two Lemmata. For the assumption in Lemma \ref{h1}, one has
 \begin{equation}
 \begin{split}
\max_{\omega \in \mathcal{O}}\frac{\|\Delta_{\omega}\mathbf{h}^n_j\|_{0}}{|\Delta \omega|} &\leq \sum^{n}_{m=1}\|\mathbf{h}^{m}_j-\mathbf{h}^{m-1}_j\|^{\gamma,\mathcal{O}}_0\frac{1}{\gamma}+ \max_{\omega \in \mathcal{O}}\frac{\|\Delta_{\omega}\mathbf{h}^0_j\|_{0}}{|\Delta\omega|}\\
& \leq \frac{1}{\gamma}\sum^{n}_{m=1}\interleave \mathbf{P}^{n} \interleave^{\gamma}_{\mathrm{s}-b,1,0} \frac{1}{j}+ |\mathrm{Op}(\langle k \rangle)|^{lip}\\
&\leq \frac{1}{\gamma \cdot j} \sum^{n}_{m=1}N^{-\mathrm{a}}_{m-1}\varepsilon + C\frac{1}{\gamma}\varepsilon \\
&\leq C
\end{split}
 \end{equation}

 For the assumption in Lemma \ref{p1}, if $N_0$ is sufficient large, one has
 \begin{equation}
 \begin{split}
 C_{\mathrm{s}}\frac{N_n^{2\tau+1}}{\gamma}\interleave \mathbf{P}^{n} \interleave^{\gamma}_{s_0,1,0} &\leq C_{\mathrm{s}} \frac{N_n^{2\tau+1}}{\gamma}N^{-\mathrm{a}}_{n-1}\interleave \mathbf{P}^{0} \interleave^{\gamma}_{\mathrm{s},1,0} \\
 & \leq \frac{1}{2}\gamma^{-1}\varepsilon \leq \frac{1}{2},
 \end{split}
 \end{equation}
 since $\frac{3}{2}(2\tau+1)-\mathrm{a}<0$.

 Now, we can apply Lemma \ref{h1} with $\mathbf{P}:=\mathbf{P}^n$ and $\mathbf{U}^+:=\mathbf{U}^{n+1}$,
 \begin{equation}
 \begin{split}
 \interleave \mathbf{U}^{n+1} \interleave^{\gamma}_{\mathrm{s}-b,1,1}&\leq \frac{N^{2\tau+1}_n}{\gamma} \interleave \mathbf{P}^{n} \interleave^{\gamma}_{\mathrm{s}-b,1,0}\\
 &\leq \interleave \mathbf{P}^{0} \interleave^{\gamma}_{s,1,0}N^{-\mathrm{a}}_{n-1}\frac{N^{2\tau+1}_n}{\gamma}.
 \end{split}
 \end{equation}
 So, the item $\mathbf{II}$ for $n+1$ is valid.

 Furthermore, from Lemma \ref{p1}, one can gets
 \begin{equation}
 \interleave \mathbf{P}^{n+1} \interleave^{\gamma}_{\mathrm{s},1,0}\leq C(\mathrm{s})\interleave \mathbf{P}^{n} \interleave^{\gamma}_{\mathrm{s},1,0}\leq C(\mathrm{s})N_{n-1}\leq N_n,
 \end{equation}
and
\begin{equation}
\begin{split}
\interleave\mathbf{P}^{n+1} \interleave^{\gamma}_{\mathrm{s}-b,1,0} &\lesssim_{\mathrm{s}} N_n^{-b}\interleave\mathbf{P}^n \interleave^{\gamma}_{\mathrm{s},1,0} +\frac{N_n^{2\tau+1}}{\gamma}\interleave\mathbf{P}^n \interleave^{\gamma}_{\mathrm{s}-b ,1,0}\interleave\mathbf{P}^n \interleave^{\gamma}_{s_0 ,1,0}\\
&\lesssim_{\mathrm{s}} N_n^{-b}N_{n-1}\interleave\mathbf{P}^0 \interleave^{\gamma}_{\mathrm{s},1,0}+\frac{N_n^{2\tau+1}}{\gamma}N^{-2\mathrm{a}}_{n-1}(\interleave\mathbf{P}^0 \interleave^{\gamma}_{\mathrm{s},1,0})^2\\
&\leq N^{-\mathrm{a}}_n \interleave\mathbf{P}^0 \interleave^{\gamma}_{\mathrm{s},1,0},
\end{split}
\end{equation}
provided
\begin{equation}\label{cs}
2C(\mathrm{s})N^{\mathrm{a}-b}_n N_{n-1} \leq 1, \quad 2C(\mathrm{s})\frac{N^{\mathrm{a}+2\tau+1}_n}{\gamma} N^{-2\mathrm{a}}_{n-1}\varepsilon \leq 1.
\end{equation}
The inequality \eqref{cs} can be verified by \eqref{as1}, by takeing $N_0$ larger enough and $\delta$ small enough. So, the item $\mathbf{IV}$ is valid for $n+1$.
 \end{proof}

\begin{cor}\label{L1}
Let $s-\tau_0-11:=\mathrm{s}>s_0+b$, and $r\in [0,\mathrm{s}-b-s_0]$. $\forall \omega \in \bigcap^{\infty}_{n=0}\mathcal{O}^{n}_{\gamma}$, the sequence of transformation
\begin{equation}
\widetilde{\Phi}_n(\theta,\omega):=\Phi_n \cdots \Phi_2 \cdot \Phi_1, \quad \Phi_n:=e^{-\mathrm{i}\mathbf{U}^n},
\end{equation}
is convergence in $ \|\cdot\|_{\mathcal{B}(\mathcal{H}^r_x \times \mathcal{H}^r_x)}^{\gamma}$ to a invertible linear operator $\widetilde{\Phi}_{\infty} $, that  fulfilling
\begin{equation}
\sup_{\theta \in \mathbb{T}^d}\|\widetilde{\Phi}^{\pm}_{\infty}(\theta)-\mathbf{Id} \|^{\gamma}_{\mathcal{B}(\mathcal{H}^r_x \times \mathcal{H}^r_{x})} \leq C\varepsilon\gamma^{-1}.
\end{equation}
\end{cor}
\begin{proof}The convergence of the transformations is a standard argument, we skip the details.
\end{proof}

\begin{cor}\label{hh1}For all $j \in \mathbb{N}$  and $\omega \in \mathcal{O}_{\gamma}$, the self-adjoint block $\{\mathbf{h}^{n}_j\}_{n\geq 0}$ is convergence in $\|\cdot\|^{\gamma}_{0}$ to a block $\mathbf{h}^{\infty}_{j}$, which is fulfils
\begin{equation}\label{h10}
\|\mathbf{h}^{\infty}_{j}-\mathbf{h}^{0}_{j}\|^{\gamma}_{0}\leq 2\varepsilon j^{-1}.
\end{equation}
\end{cor}
\begin{proof}
The convergence of the block is standard. For the bound \eqref{h1}, from \eqref{h21}, one gets
\begin{equation}
\begin{split}
\|\mathbf{h}^{\infty}_{j}-\mathbf{h}^{0}_{j}\|^{\gamma}_{0} &\leq \sum^{\infty}_{n=1}\|\mathbf{h}^n_j-\mathbf{h}^{n-1}_j\|^{\gamma}_{0} \\
&\leq \sum^{\infty}_{n=1} N^{-\mathrm{a}}_{n-2}\varepsilon j^{-1} \leq 2\varepsilon j^{-1},
\end{split}
\end{equation}
by taking $N_0$ large enough.
\end{proof}
\subsection{Measure estimate}
Set the eigenvalues of block $\mathbf{h}^{\infty}_{j}$ as $\{\lambda^{\infty}_{j,\mathfrak{a}}\}_{\mathfrak{a}\in\{1,-1\}}$, we define the set
\begin{multline} \label{non-re2}
\mathcal{O}^{\infty}_{2\gamma}:=\Big\{\omega \in \mathcal{O}_{\gamma}:|\omega \cdot \ell+\lambda^{\infty}_{i,\mathfrak{a}}+\lambda^{\infty}_{j,\mathfrak{a}'}| \geq  \frac{\gamma \langle i+ j\rangle}{\langle \ell \rangle^{\tau}}, \quad \forall (\ell,i,j)\in\mathbb{Z}^d\times \mathbb{N} \times \mathbb{N} \\
\mathfrak{a},\mathfrak{a}'\in\{1,-1\}, \quad  and \quad  |\omega \cdot \ell+\lambda^{\infty}_{i,\mathfrak{a}}-\lambda^{\infty}_{j,\mathfrak{a}'}| \geq  \frac{\gamma \langle i- j\rangle}{\langle \ell \rangle^{\tau}},\\
 \quad \forall (\ell,i,j)\in\mathbb{Z}^d\times \mathbb{N} \times \mathbb{N}, \ (\ell,i,j)\neq (0,j,j), \ \mathfrak{a},\mathfrak{a}'\in\{1,-1\}
\Big\}.
\end{multline}
\begin{lem}One has $$\mathcal{O}^{\infty}_{2\gamma}\subseteq \bigcap^{\infty}_{n=0}\mathcal{O}^{n}_{\gamma}.$$
\end{lem}
\begin{proof}It's suffice to show that  for any $n\geq0,\mathcal{O}^{\infty}_{2\gamma}\subseteq \mathcal{O}^{n}_{\gamma} $. From the definition of $\mathcal{O}^{\infty}_{2\gamma}$, one sees $\mathcal{O}^{\infty}_{2\gamma}\subseteq \mathcal{O}^{0}_{\gamma}$. For any $n>0$, from Theorem \ref{kam} and Lemma \ref{pt}, one sees that
\begin{equation}
\begin{split}
|\lambda^{\infty}_{j,\mathfrak{a}}-\lambda^{n}_{j,\mathfrak{a}}| &\leq \|\mathbf{h}^{\infty}_{j}-\mathbf{h}^{n}_{j}\|_0 \leq \sum^{\infty}_{m={n+1}}\|\mathbf{h}^{m}_{j}-\mathbf{h}^{m-1}_{j}\|_0\\
& \leq \sum^{\infty}_{m={n+1}} N^{-\mathrm{a}}_{m-2}\varepsilon j^{-1} \leq 2 N^{-\mathrm{a}}_{n-1}\varepsilon j^{-1}.
\end{split}
\end{equation}

If $\omega \in \mathcal{O}^{\infty}_{2\gamma}$ , for any $(\ell,i,j)\in\mathbb{Z}^d\times \mathbb{N} \times \mathbb{N}, \ (\ell,i,j)\neq (0,j,j), \ \mathfrak{a},\mathfrak{a}'\in\{1,-1\}$ and $|\ell| \leq N_{n-1}$, one gets
\begin{equation}\label{in1}
\begin{split}
|\omega \cdot \ell+\lambda^{n}_{i,\mathfrak{a}}-\lambda^{n}_{j,\mathfrak{a}'}| & \geq |\omega \cdot \ell+\lambda^{\infty}_{i,\mathfrak{a}}-\lambda^{\infty}_{j,\mathfrak{a}'}|-|(\lambda^{n}_{i,\mathfrak{a}}-\lambda^{\infty}_{i,\mathfrak{a}})
-(\lambda^{n}_{j,\mathfrak{a}}-\lambda^{\infty}_{j,\mathfrak{a}})|\\
&\geq \frac{2\gamma\langle i-j\rangle}{\langle\ell\rangle^{\tau}}-\frac{4\varepsilon}{jN^{\mathrm{a}}_{n-1}}\\
&\geq \frac{\gamma\langle i-j\rangle}{N_{n-1}^{\tau}}.
\end{split}
\end{equation}
The last inequality holds true, because $4\varepsilon\gamma^{-1}N^{\tau}_{n-1}\leq \langle i -j\rangle j N_{n-1}^{\mathrm{a}}$.

Also, if $\omega \in \mathcal{O}^{\infty}_{2\gamma}$ , for any $(\ell,i,j)\in\mathbb{Z}^d\times \mathbb{N} \times \mathbb{N}, \ \mathfrak{a},\mathfrak{a}'\in\{1,-1\}$ and $|\ell| \leq N_{n-1}$, one gets
\begin{equation}\label{in2}
\begin{split}
|\omega \cdot \ell+\lambda^{n}_{i,\mathfrak{a}}+\lambda^{n}_{j,\mathfrak{a}'}| & \geq |\omega \cdot \ell+\lambda^{\infty}_{i,\mathfrak{a}}+\lambda^{\infty}_{j,\mathfrak{a}'}|-|(\lambda^{n}_{i,\mathfrak{a}}-\lambda^{\infty}_{i,\mathfrak{a}})
+(\lambda^{n}_{j,\mathfrak{a}}-\lambda^{\infty}_{j,\mathfrak{a}})|\\
&\geq \frac{2\gamma\langle i+j\rangle}{\langle\ell\rangle^{\tau}}-\frac{4\varepsilon}{jN^{\mathrm{a}}_{n-1}}\\
&\geq \frac{\gamma\langle i+j\rangle}{N_{n-1}^{\tau}}.
\end{split}
\end{equation}
The last inequality holds true, because $4\varepsilon\gamma^{-1}N^{\tau}_{n-1}\leq \langle i +j\rangle j N_{n-1}^{\mathrm{a}}$. Finally, \eqref{in1} and \eqref{in2} imply that
$\mathcal{O}^{\infty}_{2\gamma}\subseteq \mathcal{O}^{n}_{\gamma}$.
\end{proof}

\begin{lem}\label{est1}Fix $\ell \in \mathbb{Z}^d\backslash \{ 0 \}$, and let $\mathcal{O}\ni \omega\mapsto \mathfrak{h}(\omega) \in \mathbb{R}$ be a Lipschitz  function fulfilling $\sup_{\omega \in \mathcal{O}}\frac{|\Delta \mathfrak{h}(\omega)|}{|\Delta \omega|}\leq \frac{1}{2}$. Define $f(\omega)=\omega\cdot \ell+\mathfrak{h}(\omega)$. Then for any $\sigma>0$. The measure of the set $\mathcal{R}:=\big\{ \omega \in \mathcal{O}| |f(\omega)|\leq \sigma \big\}$ satisfies the upper bound
\begin{equation}\label{mes1}
\mathrm{meas}(\mathcal{R}) \leq 2\frac{\sigma}{|\ell|}.
\end{equation}
\end{lem}
\begin{proof}
Fix $\ell \in \mathbb{Z}^d\backslash \{ 0 \}$, we write $\omega :=\frac{\ell}{|\ell|}\cdot r +\omega_1,\omega_1\in \mathbb{R}$ and $\omega_1 \cdot \ell =0$, then
\begin{equation}
f(\omega):=f(s)=|\ell| r+\mathfrak{h}(\omega(r)).
\end{equation}
We can obtain
$$|f(r_1-f(r_2))|\geq (|\ell|-\frac{1}{2})(r_1-r_2)\geq \frac{|\ell|}{2}(r_1-r_2),$$
such that
\begin{equation}
\mathrm{meas}\big\{r\in \mathbb{R} \big| |f(s)|\leq \sigma \big\} \leq 2\frac{\sigma}{|\ell|}.
\end{equation}
From the Fubini theorem, we can obtain \eqref{mes1}.
\end{proof}
For any $i \in \mathbb{Z}$, we known that $\langle k\rangle_{\theta,x}(i)=\mathrm{Op}(\langle k \rangle)^i_i.$ Thus, from the condition $\mathbf{II}$ and Definition \ref{defnpd}, one has
\begin{equation}
\langle k\rangle_{\theta,x}(i)(\omega)=\frac{\langle w\rangle_{\theta,x}(i)(\omega)}{\sqrt{i^2+\mathrm{m}}}=\mathfrak{c}^*(i,\omega)+\mathfrak{b}^*(i,\omega)
\end{equation}
where $\mathfrak{c}^* \in \Gamma^*:=\{\mathfrak{c}^*_{1},\cdots, \mathfrak{c}^*_{\mathrm{q}}\}$. Also, there exist an absolute positive constant $C$, such that  $|\mathfrak{b}^*(i,\omega)| \leq \frac{C}{\langle i\rangle^e}$. Take the set $\Gamma$ as $\{1,\cdots, \mathrm{q}\}$, we can define the set
\begin{multline}\label{non-re3}
\widetilde{\mathcal{O}}_{\gamma_0}:=\Big\{\omega \in \mathcal{O}:|\omega \cdot \ell+ j+\mathfrak{c}^*_{\mathbf{a}}\pm \mathfrak{c}^*_{\mathbf{a}'}|\geq \frac{\gamma_0\langle j\rangle}{\langle \ell \rangle^{\tau_0}},\quad \forall (\ell,j)\in\mathbb{Z}^{d+1}\backslash\{0\}, \ \mathbf{a}, \mathbf{a}' \in \Gamma
\Big\}.
\end{multline}

\begin{lem}\label{est2}Let $0<\gamma_0<\frac{1}{4}$, $\tau_0>d$, one has
\begin{equation}
\mathrm{meas}(\mathcal{O}\backslash \widetilde{\mathcal{O}}_{\gamma_0}) \leq C\gamma_0,
\end{equation}
where $C$ is a positive constant depending on $\mathrm{q}$.
\end{lem}
\begin{proof}If $j\neq 0$ and $\ell=0$, we known that the bound in \eqref{non-re3} hold true.

 If $j=0,\ell \neq 0 $, from Lemma \ref{est1}, the set $\mathcal{R}^{\ell,0}_{\mathbf{a},\mathbf{a}'}:=\big\{\omega \in \mathcal{O}\big||\omega\cdot\ell+ \mathfrak{c}^*_{\mathbf{a}}\pm \mathfrak{c}^*_{\mathbf{a}'}|\leq \frac{\gamma_0}{\langle\ell \rangle^{\tau_0}}\big\}$ fulfils
$$\mathrm{meas}(\mathcal{R}^{\ell,0}_{\mathbf{a},\mathbf{a}'})\leq \frac{4\gamma_0}{\langle \ell\rangle^{\tau_0+1}}. $$
Let $\mathcal{R}_1=\bigcup_{\substack{ \ell\in \mathbb{Z}^d \\ \mathbf{a}, \mathbf{a}' \in \Gamma }}\mathcal{R}^{\ell,0}_{\mathbf{a},\mathbf{a}'}$, one has
\begin{equation}
\mathrm{meas}(\mathcal{R}_1)\leq \sum_{\ell \in \mathbb{Z}^d} \sum_{\mathbf{a}, \mathbf{a}' \in \Gamma}\frac{4\gamma_0}{\langle \ell\rangle^{\tau_0+1}} \leq \sum_{\ell \in \mathbb{Z}^d}\frac{4\mathrm{p}^2\gamma_0}{\langle \ell\rangle^{\tau_0+1}}\leq C_1(\mathrm{p})\gamma_0.
\end{equation}

If $j\neq0, \ell \neq 0$ and $|j|\geq 8|\ell|$, one has
\begin{equation}
|\omega \cdot \ell+ j+\mathfrak{c}^*_{\mathbf{a}}\pm \mathfrak{c}^*_{\mathbf{a}'}|\geq |j+\mathfrak{c}^*_{\mathbf{a}}\pm \mathfrak{c}^*_{\mathbf{a}'}|-|\omega\cdot \ell|\geq \frac{1}{2}|j|-|\omega\cdot \ell|\geq\frac{1}{4}|j|\geq \frac{\gamma_0\langle j \rangle}{\langle \ell \rangle^{\tau_0+1}}
\end{equation}
Then, consider the case $1\leq |j|< 8|\ell|$. For fixed $\ell,j$, we defined the set $\mathcal{R}^{\ell,j}_{\mathbf{a},\mathbf{a}'}:=\big\{\omega \in \mathcal{O}\big||\omega\cdot\ell+ j+\mathfrak{c}^*_{\mathbf{a}}\pm \mathfrak{c}^*_{\mathbf{a}'}|\leq \frac{\gamma_0\langle j \rangle}{\langle\ell \rangle^{\tau_0}}\big\}$. Applying the Lemma \ref{est1} again, one gets
\begin{equation}
\mathrm{meas}(\mathcal{R}^{\ell,j}_{\mathbf{a},\mathbf{a}'})\leq  \frac{4\gamma_0}{\langle \ell\rangle^{\tau_0+1}}.
\end{equation}
Let $\mathcal{R}_2=\bigcap_{\substack{\ell \in \mathbb{Z}^d,|j|\leq 8|\ell| \\ \mathbf{a}, \mathbf{a}' \in \Gamma }}\mathcal{R}^{\ell,0}_{\mathbf{a},\mathbf{a}'}$, one has
\begin{equation}
\begin{split}
\mathrm{meas}(\mathcal{R}_2)&\leq \sum_{\ell \in \mathbb{Z}^d}\sum_{|j|\leq 8|\ell|} \sum_{\mathbf{a}, \mathbf{a}' \in \Gamma}\frac{4\gamma_0}{\langle \ell\rangle^{\tau_0}} \leq \sum_{\ell \in \mathbb{Z}^d}\sum_{|j|\leq 8|\ell|}\frac{4\mathrm{q}^2\gamma_0\langle j \rangle}{\langle \ell\rangle^{\tau_0+1}}\\
&\leq \sum_{\ell \in \mathbb{Z}^d}\frac{32\mathrm{m}^2\gamma_0}{\langle \ell\rangle^{\tau_0}}\leq C_2(\mathrm{q})\gamma_0
\end{split}
\end{equation}
One sees that $\mathcal{O}\backslash \widetilde{\mathcal{O}}_{\gamma_0}\subseteq \mathcal{R}_1\bigcup \mathcal{R}_2$, which finished the proof.
\end{proof}

For any $j\in\mathbb{Z}$, and $\mathfrak{a},\mathfrak{a}'\in\{1,-1\}$, we take
$$\mathrm{d}_{j,\mathfrak{a}}:=\sqrt{j^2+\mathrm{m}}+\mathfrak{c}^{*}(\mathfrak{a}j)=j+\frac{c(\mathrm{m},j)}{j}+\mathfrak{c}^{*}(\mathfrak{a}j),$$
and define the set
\begin{multline} \label{non-re2}
\widetilde{\mathcal{O}}_{\gamma_1}:=\Big\{\omega \in \widetilde{\mathcal{O}}_{\gamma_0}:|\omega \cdot \ell+\mathrm{d}_{i,\mathfrak{a}}+\mathrm{d}_{j,\mathfrak{a}'}| \geq  \frac{\gamma \langle i+ j\rangle}{\langle \ell \rangle^{\tau_1}}, \quad \forall (\ell,i,j)\in\mathbb{Z}^d\times \mathbb{N} \times \mathbb{N}   \\
\mathfrak{a},\mathfrak{a}'\in\{1,-1\}, \quad and \quad  |\omega \cdot \ell+\mathrm{d}_{i,\mathfrak{a}}-\mathrm{d}_{j,\mathfrak{a}'}| \geq  \frac{\gamma \langle i- j\rangle}{\langle \ell \rangle^{\tau_1}},\\
 \quad \forall (\ell,i,j)\in\mathbb{Z}^d\times \mathbb{N} \times \mathbb{N}, \ (\ell,i,j)\neq (0,j,j), \quad \mathfrak{a},\mathfrak{a}'\in\{1,-1\}
\Big\}.
\end{multline}

\begin{lem}\label{est3}Let $0<\gamma_1 \leq \frac{\gamma_0}{2}$ and $\tau_1>\tau_0+d$, one has
\begin{equation}\label{mea-est}
\mathrm{meas}(\widetilde{\mathcal{O}}_{\gamma_0} \backslash \widetilde{\mathcal{O}}_{\gamma_1})\leq C\frac{\gamma_1}{\gamma_0}.
\end{equation}
\end{lem}

\begin{proof}We define the set
\begin{equation}
\mathcal{U}^{\ell,i,j}=\big\{\omega \in \widetilde{\mathcal{O}}_{\gamma_0} \big| |\omega \cdot \ell+\mathrm{d}_{i,\mathfrak{a}}-\mathrm{d}_{j,\mathfrak{a}'}| \leq  \frac{\gamma_1 \langle i-j\rangle}{\langle \ell \rangle^{\tau_1}} , \forall \mathfrak{a},\mathfrak{a}' \in \{1,-1\}\big\},
\end{equation}
and
\begin{equation}
\mathcal{V}^{\ell,i,j}=\big\{\omega \in \widetilde{\mathcal{O}}_{\gamma_0} \big| |\omega \cdot \ell+\mathrm{d}_{i,\mathfrak{a}}+\mathrm{d}_{j,\mathfrak{a}'}| \leq  \frac{\gamma_1 \langle i+j\rangle}{\langle \ell \rangle^{\tau_1}} , \forall \mathfrak{a},\mathfrak{a}' \in \{1,-1\}\big\}.
\end{equation}
Let $\mathcal{U}:=\bigcup_{\substack{(\ell,i,j)\in\mathbb{Z}^d\times \mathbb{N} \times \mathbb{N}\\ (\ell,i,j) \neq (0,j,j)}}\mathcal{U}^{\ell,i,j}$ and $\mathcal{V}:=\bigcup_{(\ell,i,j)\in\mathbb{Z}^d\times \mathbb{N} \times \mathbb{N}}\mathcal{V}^{\ell,i,j}$, one  has
$$\widetilde{\mathcal{O}}_{\gamma_0}\backslash \widetilde{\mathcal{O}}_{\gamma_1} \subseteq \mathcal{U}\bigcup \mathcal{V}.$$
We consider the measure estimate of set $\mathcal{P}$, estimateing the measure  of set $\mathcal{Q}$ is relatively simple.

\textbf{Case 1}: If $\ell=0$ and $i\neq j$, one has
\begin{equation}
|\mathrm{d}_{i,\mathfrak{a}}-\mathrm{d}_{j,\mathfrak{a}'}| \geq \frac{1}{2}|i-j| \geq \gamma_1 \langle i-j\rangle.
\end{equation}

\textbf{Case 2}: If $\ell \neq 0$ and $i= j$, one has
\begin{equation}
|\omega\cdot\ell+\mathfrak{c}^*(\mathfrak{a}j)-\mathfrak{c}^*(\mathfrak{a}'j)| \geq \frac{\gamma_0}{\langle\ell \rangle^{\tau_0}} \geq \frac{\gamma_1}{\langle\ell \rangle^{\tau_1}}
\end{equation}

\textbf{Case 3}: If $\ell\neq 0,i\neq j$ and  $|i-j|> 8|\ell|$, one can obtain
\begin{equation}
|\omega\cdot\ell+\mathrm{d}_{i,\mathfrak{a}}-\mathrm{d}_{j,\mathfrak{a}'}| \geq \frac{1}{2}|i-j|-|\omega\cdot \ell | \geq \frac{1}{4}|i-j| \geq \frac{\gamma_1 \langle i-j\rangle}{\langle \ell\rangle^{\tau_1}}
\end{equation}

\textbf{Case 4}: Let $|i-j|\leq 8|\ell|$ and $i< j$, we assume that
\begin{equation}
\langle i \rangle \langle i-j\rangle \geq \frac{4\mathrm{m}\langle\ell \rangle^{\tau_0}}{\gamma_0},
\end{equation}
then
\begin{equation}
\begin{split}
|\omega \cdot \ell+\lambda^{\infty}_{i,\mathfrak{a}}+\lambda^{\infty}_{j,\mathfrak{a}'}|&\geq |\omega\cdot\ell+i-j+\mathfrak{c}^{*}_{i,\mathfrak{a}}-\mathfrak{c}^{*}_{j,\mathfrak{a}'}|-\frac{2\mathrm{m}}{\langle i \rangle}\\
&\geq \frac{\gamma_0\langle i -j\rangle}{\langle\ell \rangle^{\tau_0}}-\frac{2m}{\langle i \rangle}\\
&\geq \frac{\gamma_0\langle i -j\rangle}{2\langle\ell \rangle^{\tau_0}}.
\end{split}
\end{equation}
Therefore, we restrict ourself to the case $i<j$ and $\langle i \rangle\langle i-j\rangle \leq \frac{4\mathrm{m}\langle\ell \rangle^{\tau_0}}{\gamma_0}$. The same arguments can be extended to the symmetric case $j> i$ and $\langle j \rangle\langle i-j\rangle \leq \frac{4\mathrm{m}\langle\ell \rangle^{\tau_0}}{\gamma_0}$.

From Lemma \ref{est1}, we known that for any $\ell\neq 0$
\begin{equation}\label{u11}
\mathrm{meas}(\mathcal{U}^{\ell,i,j})\leq \frac{8\gamma_1\langle i-j \rangle}{\langle\ell \rangle^{\tau_1}}
\end{equation}

Now, we define the index set of $(\ell,i,j)$, that is
$$\mathcal{E}:=\{|i-j|\leq 8|\ell| \}\bigcap \Big(\{i\leq j,\langle i \rangle \langle i-j\rangle \leq \frac{4\mathrm{m}\langle\ell \rangle^{\tau_0}}{\gamma_0}\}\bigcup\{j\leq i, \langle j \rangle \langle i-j\rangle \leq \frac{4\mathrm{m}\langle\ell \rangle^{\tau_1}}{\gamma_0}\}\Big).$$
Since $\mathcal{U}=\bigcup_{(\ell,i,j)\in \mathcal{E}}\mathcal{U}^{\ell,i,j} $, from \eqref{u11}, one gets
\begin{equation}
\begin{split}
\mathrm{meas}(\mathcal{U})&\leq \sum_{(\ell,i,j)\in \mathcal{E} }\mathrm{meas}(\mathcal{U}^{\ell,i,j})\\
&\leq 16\gamma_1 \sum_{\ell \neq 0}\sum_{\substack{i< j \\ \langle i \rangle\langle i-j\rangle \leq \frac{4\mathrm{m}\langle\ell \rangle^{\tau_0}}{\gamma_0}}} \sum_{|i-j|\leq 8|\ell|}
\frac{\langle i -j \rangle}{\langle\ell \rangle^{\tau_1+1}}\\
&\leq 16\gamma_1 \sum_{\ell \neq 0}\sum_{\substack{j-i:=k \\ k\leq 8|\ell|}} \sum_{\langle i \rangle \leq \frac{4\mathrm{m}\langle\ell \rangle^{\tau_0}}{k \gamma_0}} \frac{k}{\langle\ell \rangle^{\tau_1=1}}\\
&\leq 64 \mathrm{m}\frac{\gamma_1}{\gamma_0}\sum_{\ell \neq 0}\sum_{\substack{j-i:=k \\ k\leq 8|\ell|}} \frac{1}{\langle\ell \rangle^{\tau-\tau_0+1}}\\
&\lesssim 512\mathrm{m} \frac{\gamma_1}{\gamma_0}\sum_{\ell \neq 0}\frac{1}{\langle \ell \rangle^{\tau-\tau_0}} \\
&\lesssim C\frac{\gamma_1}{\gamma_0},
\end{split}
\end{equation}
provided $\tau>d+\tau_0.$ The same computation hold for the set $\mathcal{V}$. Hence, we conclude the estimate \eqref{mea-est}.
\end{proof}

From the Lemma \ref{pt} and Corollary \ref{hh1}, for any $j \in \mathbb{N}$ and $\mathfrak{a} \in \{1.-1\}$,  the final eigenvalues $\lambda^{\infty}_{j,\mathfrak{a}}$ fulfils
\begin{equation}
\begin{split}
\lambda^{\infty}_{j,\mathfrak{a}}:&=\lambda^{\infty}_{j,\mathfrak{a}}+\varepsilon^{\infty}_{j,\mathfrak{a}}(\omega)\\
&=\sqrt{j^2+\mathrm{m}^2}+\langle k \rangle_{\theta,x}(\mathfrak{a}j)+\varepsilon^{\infty}_{j,\mathfrak{a}}(\omega)\\
&=j+\frac{c(\mathrm{m},j)}{j}+\mathfrak{c}^{*}(\mathfrak{a}j)+\mathfrak{b}^{*}(\mathfrak{a}j)+\varepsilon^{\infty}_{j,\mathfrak{a}}(\omega),
\end{split}
\end{equation}
where $$|\mathfrak{b}^{*}(\mathfrak{a}j)|^{\gamma}\leq \frac{c\epsilon}{\langle j\rangle^e}, \ |\varepsilon^{\infty}_{j,\mathfrak{a}}(\omega)|^{\gamma}\leq \frac{2\varepsilon}{\langle j\rangle}.$$
Take $\rho:=\min\{1,e\}$, one gets

\begin{equation}
|\mathfrak{b}^{*}(\mathfrak{a}j)+\varepsilon^{\infty}_{j,\mathfrak{a}}(\omega)|^{\gamma}\leq \frac{c_1\varepsilon}{\langle j \rangle^{\rho}}, \ c_1=2+c.
\end{equation}
\begin{lem}\label{est4}Let $0<\gamma<\frac{\gamma_1}{2}$ and $\tau> \max \{d+\frac{\tau_1}{\rho}-\frac{1}{\rho}, d+\frac{\tau_0}{\rho}-1\} $, one has that
\begin{equation}\label{mea-est2}
\mathrm{meas}(\widetilde{\mathcal{O}}_{\gamma_1} \backslash \mathcal{O}^{\infty}_{2\gamma}) \leq C \gamma.
\end{equation}
\end{lem}
\begin{proof}We define the set
\begin{equation}
\mathcal{P}^{\ell,i,j}=\big\{\omega \in \widetilde{\mathcal{O}}_{\gamma_0} \big| |\omega \cdot \ell+\lambda^{\infty}_{i,\mathfrak{a}}-\lambda^{\infty}_{j,\mathfrak{a}'}| \leq  \frac{2\gamma \langle i-j\rangle}{\langle \ell \rangle^{\tau}} , \forall \mathfrak{a},\mathfrak{a}' \in \{1,-1\}\big\},
\end{equation}
and
\begin{equation}
\mathcal{Q}^{\ell,i,j}=\big\{\omega \in \widetilde{\mathcal{O}}_{\gamma_0} \big| |\omega \cdot \ell+\lambda^{\infty}_{i,\mathfrak{a}}+\lambda^{\infty}_{j,\mathfrak{a}'}| \leq  \frac{2\gamma \langle i+j\rangle}{\langle \ell \rangle^{\tau}} , \forall \mathfrak{a},\mathfrak{a}' \in \{1,-1\}\big\}.
\end{equation}
Let $\mathcal{P}:=\bigcup_{\substack{(\ell,i,j)\in\mathbb{Z}^d\times \mathbb{N} \times \mathbb{N} \\ (\ell,i,j) \neq (0,j,j)}}\mathcal{P}^{\ell,i,j}$ and $\mathcal{Q}:=\bigcup_{(\ell,i,j)\in\mathbb{Z}^d\times \mathbb{N} \times \mathbb{N}}\mathcal{Q}^{\ell,i,j}$, one  gets
$$\widetilde{\mathcal{O}}_{\gamma_0}\backslash \mathcal{O}^{\infty}_{2\gamma} \subseteq \mathcal{P}\bigcup \mathcal{Q}.$$
We focus on the measure estimate of set $\mathcal{P}$, it's relatively simple to estimate the measure of set $\mathcal{Q}$.

\textbf{case 1}: If $\ell=0$ and $i\neq j$, one has
\begin{equation}
|\lambda^{\infty}_{i,\mathfrak{a}}+\lambda^{\infty}_{j,\mathfrak{a}'}| \geq \frac{1}{2}|i-j| \geq 2\gamma \langle i-j\rangle.
\end{equation}

\textbf{case 2}: If $\ell\neq 0$ and $i=j$, let $\langle j \rangle^{\rho}>c(\varepsilon, \gamma_0)\langle \ell \rangle^{\tau_0}$, one has
\begin{equation}\label{p11}
\begin{split}
|\omega\cdot\ell+\lambda^{\infty}_{j,\mathfrak{a}}+\lambda^{\infty}_{j,\mathfrak{a}'}| & \geq |\omega \cdot \ell+\mathfrak{c}^{*}(\mathfrak{a}j)-\mathfrak{c}^{*}(\mathfrak{a}'j)|-\frac{2c_1 \varepsilon}{\langle j \rangle^{\rho}}\\
&\geq \frac{\gamma_0}{\langle \ell \rangle^{\tau_0}}-\frac{2c_1 \varepsilon}{\langle j \rangle^{\rho}}\\
&\geq \frac{\gamma_0}{2 \langle \ell \rangle^{\tau_0}}.
\end{split}
\end{equation}

Let $\mathcal{P}_{1}=\bigcup_{\substack{\ell \in \mathbb{Z}^d,j\in \mathbb{N}\\ (\ell,j,j) \neq (0,j,j)}}\mathcal{P}^{\ell,j,j}$, from \eqref{p1}, one has
$\mathcal{P}_{1}=\bigcup_{\substack{\ell \in \mathbb{Z}^d,\langle j \rangle^{\rho} \leq  c(\varepsilon,\gamma_0)\langle \ell \rangle^{\tau_0}\\ (\ell,j,j) \neq (0,j,j)}}\mathcal{P}^{\ell,j,j}.$

From Lemma \ref{est1}, for any $\ell \neq 0$, one gets
\begin{equation}
\mathrm{meas}(\mathcal{P}^{\ell,j,j}) \leq \frac{16\gamma}{\langle \ell \rangle^{\tau+1}}.
\end{equation}
Then
\begin{equation}\label{p1.1}
\begin{split}
\mathrm{meas}(\mathcal{P}_1)&\leq \sum_{\ell \in \mathbb{Z}^d \backslash \{0\}} \sum_{\langle j \rangle^{\rho} \leq  c(\varepsilon,\gamma_0)\langle \ell \rangle^{\tau_0}} \frac{16\gamma}{\langle \ell \rangle^{\tau+1}}\\
&\leq 16\gamma(c(\varepsilon,\gamma_0))^{\frac{1}{\rho}}\sum_{\ell \in \mathbb{Z}^d \backslash \{0\}}\frac{1}{\langle \ell \rangle^{\tau-\frac{\tau_0}{\rho}+1}}\\
&\leq \tilde{c}(\varepsilon,\gamma_0)\gamma,
\end{split}
\end{equation}
provided $\tau-\frac{\tau_0}{\rho}+1>d$.

\textbf{case 3}: If $\ell=0$ and $i\neq j$. Let $|i-j|\geq 8|\ell|$, we can get
\begin{equation}
|\omega\cdot\ell+\lambda^{\infty}_{i,\mathfrak{a}}+\lambda^{\infty}_{j,\mathfrak{a}'}| \geq \frac{1}{2}|i-j|-|\omega\cdot \ell | \geq \frac{1}{4}|i-j| \geq \frac{2\gamma \langle i-j\rangle}{\langle \ell\rangle^{\tau}}
\end{equation}

\textbf{case 4}:Let $|i-j|\leq 8|\ell|$ and $i< j$, we assume that
\begin{equation}
\langle i \rangle^{\rho}\langle i-j\rangle \geq c(\varepsilon,\gamma_1)\langle\ell \rangle^{\tau_1}£¬
\end{equation}
then
\begin{equation}
\begin{split}
|\omega \cdot \ell+\lambda^{\infty}_{i,\mathfrak{a}}+\lambda^{\infty}_{j,\mathfrak{a}'}|&\geq |\omega\cdot\ell+\mathrm{d}_{i,\mathfrak{a}}-\mathrm{d}_{j,\mathfrak{a}'}|-\frac{2c_1\varepsilon}{\langle i \rangle^{\rho}}\\
&\geq \frac{\gamma_1\langle i -j\rangle}{\langle\ell \rangle^{\tau_1}}-\frac{2c_1\varepsilon}{\langle i \rangle^{\rho}}\\
&\geq \frac{\gamma_1\langle i -j\rangle}{2\langle\ell \rangle^{\tau_1}}.
\end{split}
\end{equation}
Therefore, we restrict ourself to the case $i< j$ and $\langle i \rangle^{\rho}\langle i-j\rangle \leq c(\varepsilon,\gamma_1)\langle\ell \rangle^{\tau_1}$. The same arguments can be extended to the symmetric case $j< i$ and $\langle j \rangle^{\rho}\langle i-j\rangle \leq c(\varepsilon,\gamma_1)\langle\ell \rangle^{\tau_1}$.
From Lemma \ref{est1}, we known that for any $\ell\neq 0$ and $i\neq j$
\begin{equation}\label{u1}
\mathrm{meas}(\mathcal{P}^{\ell,i,j})\leq \frac{16\gamma \langle i-j \rangle}{\langle\ell \rangle^{\tau+1}}.
\end{equation}

Now, we define the index set of $(\ell,i,j)$, that is
$$\mathcal{T}:=\{|i-j|\leq 8|\ell| \}\bigcap \Big(\{i<j,\langle i \rangle^{\rho}\langle i-j\rangle \leq c(\varepsilon,\gamma_1)\langle\ell \rangle^{\tau_1}\}\bigcup\{j< i, \langle j \rangle^{\rho}\langle i-j\rangle \leq c(\varepsilon,\gamma_1)\langle\ell \rangle^{\tau_1}\}\Big).$$
Let $\mathcal{U}_2=\bigcup_{(\ell,i,j) \in \mathcal{T}}$ , one gets
\begin{equation}\label{p1.2}
\begin{split}
meas(\mathcal{P}_2)&\leq \sum_{(\ell,i,j)\in \mathcal{T} }meas(\mathcal{P}^{\ell,i,j})\\
&\leq 16\gamma \sum_{\ell \neq 0}\sum_{\substack{i< j \\ \langle i \rangle^{\rho}\langle i-j\rangle \leq c(\varepsilon,\gamma_1)\langle\ell \rangle^{\tau_1}}} \sum_{|i-j|\leq 8|\ell|}
\frac{\langle i -j \rangle}{\langle\ell \rangle^{\tau+1}}\\
&\leq 16 \gamma \sum_{\ell \neq 0} \sum_{\substack{j-i:=k \\ k\leq 8|\ell|}}\sum_{\langle i\rangle^{\rho} \leq \frac{c(\varepsilon,\gamma_1)\langle \ell \rangle^{\tau_1}}{k}}\frac{k}{\langle\ell \rangle^{\tau+1}}\\
&\leq 16\gamma (c(\varepsilon,\gamma_1))^{\frac{1}{\rho}} \sum_{\ell \neq 0} \sum_{\substack{j-i:=k \\ k\leq 8|\ell|}}\frac{k^{1-\frac{1}{\rho}}}{\langle\ell \rangle^{\tau-\frac{\tau_1}{\rho}+1}}\\
& \leq 16\cdot 8^{1-\frac{1}{\rho}} (c(\varepsilon,\gamma_1))^{\frac{1}{\rho}} \gamma  \sum_{\ell \neq 0} \frac{1}{\langle\ell \rangle^{\tau-\frac{\tau_1}{\rho}+\frac{1}{\rho}}}\\
&\leq  \tilde{c}(\varepsilon,\gamma_1) \gamma,
\end{split}
\end{equation}
provided $\tau-\frac{\tau_1}{\rho}+\frac{1}{\rho}>d$. The bounds \eqref{p1.1} and \eqref{p1.2} imply that
$$ \mathrm{meas}(\mathcal{P}) \leq \tilde{c}(\varepsilon,\gamma_0,\gamma_1)\gamma. $$
The same computation hold for the set $\mathcal{Q}$. Hence, we conclude the  measure estimate \eqref{mea-est2}.
\end{proof}
\begin{prop}\label{mea}One has
\begin{equation}
\mathrm{meas}(\mathcal{O} \backslash \mathcal{O}^{\infty}_{2\gamma}) \leq \tilde{C}\gamma^{\frac{1}{3}}.
\end{equation}
\end{prop}
\begin{proof}
From the definitions of sets $\mathcal{O}_{\gamma},\widetilde{\mathcal{O}}_{\gamma_0},\widetilde{\mathcal{O}}_{\gamma_1},\mathcal{O}^{\infty}_{2\gamma}$, one gets
\begin{equation}
\mathcal{O} \backslash \mathcal{O}^{\infty}_{2\gamma}=(\mathcal{O}\backslash \mathcal{O}_{\gamma})\bigcup(\mathcal{O}_{\gamma}\backslash \widetilde{\mathcal{O}}_{\gamma_0})\bigcup(\widetilde{\mathcal{O}}_{\gamma_0}\backslash\widetilde{ \mathcal{O}}_{\gamma_1})\bigcup (\widetilde{\mathcal{O}}_{\gamma_1}\backslash \mathcal{O}^{\infty}_{2\gamma}).
\end{equation}
Let $\gamma_0=\gamma^{\frac{1}{3}}$ and  $\gamma_1=\gamma^{\frac{2}{3}}$, from Lemmata \ref{est2}, \ref{est3}, \ref{est4},  we can get
\begin{equation}
\mathrm{meas}(\mathcal{O} \backslash \mathcal{O}^{\infty}_{2\gamma}) \leq C_0\gamma+C_1\gamma^{\frac{1}{3}}+C_2\gamma^{\frac{1}{3}}+C_3 \gamma
\leq \tilde{C}\gamma^{\frac{1}{3}},
\end{equation}
where $\tilde{C}:=4 \cdot \max\{C_0,C_1,C_2,C_3\}$.
\end{proof}
\section {Proof of main Theorem \ref{main1} and Corollary \ref{main2}.}
We define the composition operator
\begin{equation}
\mathbf{A}(\theta,\omega):=  \widetilde{\Phi}_{\infty}(\theta,\omega)\circ \mathbf{V}(\theta,\omega)
\end{equation}
where $\mathbf{V}(\theta,\omega)$ is defined in Remark \ref{regu} and $ \widetilde{\Phi}_{\infty}(\theta,\omega)$ is defined in Corollary \ref{L1}. We also define the constants
$$\bar{s}:=s_0+\tau_0+11+\mathrm{b}$$
and for any $s>\bar{s}$, we define
$$\mathfrak{R}_s:=s-\tau_0-11-\mathrm{b}-s_0,$$
where we recall the definitions  in \eqref{as1}. From Lemmata \ref{link2}, \ref{MS2}, \ref{G1} and Theorem \ref{kam}, one gets that for $\varepsilon \gamma^{-1}\leq \delta_{s}$, for any $\theta\in \mathbb{T}^d$ and $\omega \in \mathcal{O}^{\infty}_{2\gamma}$,the maps $\mathbf{A}^{\pm}(\theta,\omega)$ are bounded and invertible with
\begin{equation}\label{A}
\mathbf{A}^{\pm}(\theta,\omega):(\mathcal{H}^{r}_x\times\mathcal{H}^{r}_x )\mapsto (\mathcal{H}^{r}_x\times\mathcal{H}^{r}_x ),
\end{equation}
for any $0\leq r \leq \mathfrak{R}_s$.

Also, for any $\omega \in \mathcal{O}^{2\gamma}_{\infty}$, by the change of variables $\mathbf{q}:=\mathbf{A}(\omega t) \mathbf{v}$, the Cauchy problem
$$\left\{
\begin{aligned}
&\mathrm{i}\partial_{t}\mathbf{q}(t)=\mathbf{H}(t)\mathbf{q}(t), \\
&\mathbf{q}(0,x)=(q(0,x),\bar{q}(0,x)),
\end{aligned}
\right.
$$
is transformed into
$$\left\{
\begin{aligned}
&\mathrm{i}\partial_{t}\mathbf{v}(t)=\mathbf{H}^{\infty}_0\mathbf{v}(t) \\
&\mathbf{v}(0,x)=(v(0,x),\bar{v}(0,x))
\end{aligned}
\right. ,
\quad \mathbf{v}(0,x)=\mathbf{A}^{-1}(0,\omega)\mathbf{q}(0,x), $$
where the operator $\mathbf{H}^{\infty}_0=\left(
\begin{array}{cc}
\mathcal{H}^{\infty}_0 & 0 \\
0 & -\overline{\mathcal{H}}^{\infty}_0
\end{array}
\right )$ is defined in Corollary \ref{hh1}. Then, we can consider the Cauchy problem
$$\left\{
\begin{aligned}
&\mathrm{i}\partial_{t}v(t)=\mathcal{H}^{\infty}_0v(t), \\
&v_0(x)=v(0,x).
\end{aligned}
\right.
$$
Since the operator $\mathcal{H}^{\infty}_0 $ is block-diagonal and self-adjoint, we can verified that
\begin{equation}
\partial_{t}\|v(t,x)\|^2_{\mathcal{H}^r_x}=-(\mathrm{i}(\mathcal{H}^{\infty}_0-(\mathcal{H}^{\infty}_0)^*)\langle D\rangle^rv,\langle D\rangle^rv)=0,
\end{equation}
which implies that
\begin{equation}\label{v1}
\|v(t,x)\|_{\mathcal{H}^r_x}=\|v(0,x)\|_{\mathcal{H}^r_x}.
\end{equation}
By \eqref{A} and \eqref{v1}, we can get
\begin{equation}
\|\mathbf{q}(0,x)\|_{\mathcal{H}^r_x \times \mathcal{H}^r_x} \lesssim_r\|\mathbf{q}(t,x)\|_{\mathcal{H}^r_x \times \mathcal{H}^r_x}  \lesssim_r \|\mathbf{q}(0,x)\|_{\mathcal{H}^r_x \times \mathcal{H}^r_x}.
\end{equation}
Set $\gamma=\varepsilon^{a},0<a<1$ and $\mathcal{O}_{\epsilon}=\mathcal{O}_{2\gamma}^{\infty}$, the Proposition \ref{mea} implies that $$\lim_{\epsilon\rightarrow 0}\mathrm{meas}(\mathcal{O}\backslash \mathcal{O}_{\epsilon} )=0.$$
\section{Appendix}

\subsection{Properties of self-adjoint matrix}\

In this section, we recall some well known facts about self-adjoint  operator in the finite dimension Hilbert space $\mathcal{H}$. Let $\mathcal{H}$ be a
finite dimensional Hilbert space of dimension $\mathfrak{n}$ equipped by the inner product
$(, )_\mathcal{H}$. For any self-adjoint  operator $A$, we order its eigenvalues as
$spec(A) := {\lambda_1(A) \leq \lambda_2(A) \leq \cdots \leq \lambda_\mathfrak{n}(A)}.$

\begin{prop}\label{pt}(Weyl's Perturbation Theorem)(\cite{R97}, Theorem III.2.1) Let $A$ and $B$ be self-adjoint matrices. Then
\begin{equation}
|\lambda_k(A)-\lambda_k(B)|\leq \|A-B\|_{0},\  \forall k \in 1,\cdots,\mathfrak{n}.
\end{equation}
\end{prop}
\begin{prop}\label{pt2}(\cite{R97}, Theorem VII.2.8) Let $A$ and $B$ be self-adjoint matrices, and let $\delta=dist(\sigma(A),\sigma(B))$. Then the solution $X$ of the equation $AX-XB=Y$ satisfies the inequality
\begin{equation}
\|X\|_{0}  \leq \frac{C}{\delta} \|Y\|_{0}.
\end{equation}
\end{prop}

\section*{Acknowledgements}
The work is supported by the Jiangsu Province Postdoctoral Foundation(No.2021K163B), China Postdoctoral Foundation (No.2021M692717), and the National Natural Science Foundation of China
(No.12101542).

\bibliographystyle{abbrv} 

\begin{thebibliography}{10}
\bibitem{bal19}
P.~Baldi, M.~Berti, E.~Haus, R.~Montalto.
\newblock Time quasi-periodic gravity water waves in finite depth.
\newblock {\em Invent. Math.}, 214, 739-911, 2018.

\bibitem{Baldi2}
P.~Baldi, M.~Berti, E.~Haus, R.~Montalto.
\newblock  KAM for quasi-linear and fully nonlinear forced perturbations of Airy equation.
\newblock  {\em Mathematische. Annalen.}, 359, 471-536, 2014.


\bibitem{Bam01}
D.~Bambusi and S.~Graffi.
\newblock Time Quasi-periodic unbounded perturbations of Schr\"odinger operators and KAM methods.
\newblock {\em Comm. Math. Phys.}, 219: 465-480, 2001.
\bibitem{Bam18}
D.~Bambusi.
\newblock Reducibility of 1-d Schr\"odinger equation with time quasiperiodic unbounded perturbation,I.
\newblock{\em Trans. Amer. Math. Soc.}, 370: 1823-1865, 2018.

\bibitem{Bam171}
D.~Bambusi.
\newblock Reducibility of 1-d Schr\"odinger equation with time quasiperiodic unbounded perturbation, II.
\newblock{\em Comm. Math. Phys.}, 353: 353-378, 2017.

\bibitem{Bam19}
D.~Bambusi and R.~Montalto.
\newblock Reducibility of 1-d Schr\"odinger equation with time quasiperiodic unbounded perturbation, III.
\newblock{\em  J. Math. Phys.}, 59, 2018.



\bibitem{Bam018}
D. Bambusi, B. Gr\'ebert, A. Maspero, and D. Robert.
\newblock Reducibility of the quantum harmonic oscillator in d-dimensions with polynomial time-dependent perturbation.
\newblock {\em Anal. PDE.,} 11: 775-799, 2018.

\bibitem{Bam20}
D.~Bambusi, B.~Langella  and R.~Montalto.
\newblock Growth of Sobolev norms for unbounded perturbations of the Laplacian on flat tori.
\newblock Preprint, arXiv:2012.02654.

\bibitem{Bam22}
D.~Bambusi, B.~Langella  and R.~Montalto.
\newblock Spectral asymptotics of all the eigenvalues of Schr\"odinger operators on flat tori.
\newblock {\em Nonlinear. Anal.}, 216: 2022.

\bibitem{Bam21}
D.~Bambusi, B.~Gr\'ebert, A.~Maspero and D.~Robert.
\newblock Growth of Sobolev norms for abstract linear Schr\"odinger Equations.
\newblock {\em  J. Eur. Math. Soc.}, 23:557-583, 2021.

\bibitem{Berti13}
M.~Berti, L.~Biasco, M.~Procesi.
\newblock KAM theory for the Hamiltonian derivative wave equation.
\newblock {\em Ann. Sci. \'Ec. Norm. Sup\'er.}, 46; 301-373, 2013.

\bibitem{berti}
M.~Berti, L.~Corsi, M.~Procesi.
\newblock  An abstract Nash Moser theorem and quasi-periodic solutions for NLW and NLS on compact Lie groups and Homogeneous manifolds.
\newblock {\em Comm. Math. Phys.}, 334: 1413-1454, 2015.

\bibitem{berti19}
M.~Berti, A.~Maspero.
\newblock Long time dynamics of Schr\"odinger and wave equations on flat tori
\newblock {\em J. Differential. Equations.}, 267: 1167-1200, 2019.

\bibitem{bert20}
M.~Berti, R.~Montalto.
\newblock Quasi-Periodic Standing Wave Solutions of Gravity-Capillary Water Waves.
\newblock {\em Mem. Amer. Math. Soc.}, 263 : 2020.

\bibitem{R97}
R.~Bhatia.
\newblock Matrix Analysis.
\newblock{\em Springer-Verlag New York}, 1997.

\bibitem{Bou99}
J.~Bourgain.
\newblock On growth of Sobolev norms in linear Schr\"odinger equations with smooth time dependent potential.
\newblock {\em J. Anal. Math.}, 77: 315-348, 1999.

\bibitem{Bou99.1}
J.~Bourgain.
\newblock Growth of Sobolev norms in linear Schr\"odinger equations with quasi-periodic potential.
\newblock {\em Comm. Math. Phys.}, 204:207-247, 1999.


\bibitem{Delort2010}
\newblock J.~ Delort.
\newblock Growth of Sobolev norms of solutions of linear Schr\"odinger equations on some compact manifolds.
\newblock {\em Int. Math. Res. Not.}, 12: 2305-2328,2010.


\bibitem{L.H 09}
L.~Eliasson and S.~Kuksin.
\newblock  On reducibility of Schr\"odinger equation with quasi periodic in time potential.
\newblock {\em Comm. Math. Phys.}, 286: 125-135, 2009.


\bibitem{Fang14}
D.~Fang,  Z.~Han, W.-M.~Wang.
\newblock Bounded Sobolev norms for Klein-Gordon equations under non-resonant perturbation.
\newblock {\em J. Math. Phys.}, 55: 2014.

\bibitem{L.M2019}
L.~Franzoi and A.~Maspero.
\newblock Reducibility for a fast-driven linear Klein-Gordon equation.

\newblock{\em Annali. di. Matematica. Pura. ed. Applicata.}, 198: 1407-1439, 2019.




\bibitem{F.G2019}
R.~Feola, B.~Gr\'ebert.
\newblock Reducibility of Schr\"odinger equation on sphere.
\newblock {\em Int. Math. Res. Notices.}, 19: 15082-15120, 2021.


\bibitem{F.G2020}
R.~Feola, B.~Gr\'ebert, T.~Nguyen.
\newblock Reducibility of Schr\"odinger equation on a Zoll manifold with unbounded potential.
\newblock {\em J. Math. Phys.}£¬ 61£» 2020.


\bibitem{B12}
B.~Gr\'ebert and  L.~Thomann.
\newblock KAM for the Quantum Harmonic Oscillator.
\newblock {\em Commun. Math. Phys.}, 307: 383-427, 2011.

\bibitem{B24}
B.~Gr\'ebert and E.~Paturel.
\newblock On reducibility of Quantum Harmonic oscillator on $\mathbb{R}^d$ with quasiperiodic in time potential.
\newblock {\em  Ann. Fac. Sci. Toulouse.}, 6:977-1014, 2019.

\bibitem{Li2020}
J.~Li.
\newblock Reducibility, Lyapunov exponent, pure point spectra property for quasi-periodic wave operator.
\newblock {\em Taiwanese J. Math.}£¬ 24; 377-411, 2020.

\bibitem{Liu09}
J.~Liu and X.~Yuan.
\newblock Spectrum for quantum Duffing oscillator and small-divisor equation with large-variable coefficient.
\newblock {\em Comm. Pure. Appl. Math.}, 63: 1145-1172, 2010.

\bibitem{Mas2017}
A.~Maspero and  D.~Robert.
\newblock On time dependent Schr\"odinger equations: global well-posedness and growth of Sobolev norms.
\newblock {\em J. Funct. Anal.}, 273£»721-781, 2017.

\bibitem{Mon2017}
R.~Montalto.
\newblock Quasi-periodic solutions of forced Kirchhoff equation.
\newblock{\em NoDEA Nonlinear. Differential. Equations. Appl.}, 24: 2017.

\bibitem{Mon2018.1}
R.~Montalto.
\newblock On the growth of Sobolev norms for a class of linear Schr\"odinger equations on the torus with superlinear dispersion.
\newblock {\em Asymptot. Anal.}, 108: 85-114, 2018.

\bibitem{Mo2019}
R.~Montalto.
\newblock A reducibility result for a class of linear  wave equation  on $\T^d$.
\newblock {\em Int. Math. Res. Notices.}, 6: 1788-1862, 2019.




\bibitem{Mo20191}
R.~Montalto.
\newblock Growth of Sobolev norms for time dependent periodic Schr\"odinger equations with sublinear dispersion.
\newblock{\em J. Differential. Equations.}, 266: 4953-4996, 2019.


\bibitem{S19}
Y.~Sun, J.~Li and  B.~Xie.
\newblock Reducibility for wave equations of finitely smooth potential with periodic boundary conditions.
\newblock{\em J. Differential. Equations.}, 266: 2762-2804, 2019.

\bibitem{S21}
Y.~Sun, J.~Li.
\newblock Reducibility of relativistic Schr\"odinger equation with unbounded perturbations.
\newblock {\em J. Differential. Equations.}, 286:215-247, 2021.

\bibitem{Liang2018}
Z.~Liang, X.~Wang.
\newblock On reducibility of 1d wave equation with quasiperiodic in time potentials.
\newblock {\em J. Dynam. Differential. Equations.}, 30: 957-978, 2018.

\bibitem{Wang08}
W.-M.~Wang.
\newblock Logarithmic bounds on Sobolev norms for time dependent linear Schr\"odinger equations.
\newblock {\em Comm. Partial. Differential. Equations.}, 33: 2164-2179,2008.
\end{thebibliography}

\end{document}